\documentclass{amsart}
\usepackage{color}
\usepackage[colorlinks,linkcolor=blue,citecolor=blue,urlcolor=blue]{hyperref}
\usepackage{mathrsfs}
\usepackage{enumerate}
\usepackage{graphicx}

\allowdisplaybreaks

\newcommand{\eq}{\begin{equation}}
\newcommand{\en}{\end{equation}}

\makeatletter

\newcommand{\Rmnum}[1]{\expandafter\@slowromancap\romannumeral #1@}
\makeatother

\newtheorem{thm}{Theorem}[section]
\newtheorem{prop}[thm]{Proposition}
\newtheorem{lem}[thm]{Lemma}
\newtheorem{cor}[thm]{Corollary}
\theoremstyle{definition}
\newtheorem{remark}[thm]{Remark}
\newtheorem{defn}[thm]{Definition}

\newtheorem{example1}[thm]{Example}
\newtheorem*{hypA}{Hypothesis \textbf{A}}
\newtheorem*{hypB}{Hypothesis \textbf{B}}
\newtheorem*{assB}{Assumption \textbf{B}}

\numberwithin{equation}{section}

\newcommand{\reals}{\mathbb{R}}

\newcommand{\nats}{\mathbb{N}}

\newcommand{\ra}{\rightarrow}
\renewcommand{\P}{\mathbb{P}}

\newcommand{\E}{\mathbb{E}}

\newcommand{\eps}{\varepsilon}
\renewcommand{\and}{ \quad \text{ and } \quad }

\begin{document}

\title[Shocks in Burgers turbulence with L\'evy noise
  initial data]{Structure of shocks in Burgers turbulence \\
with L\'evy noise initial data}

\author{Joshua Abramson}
\address{J.~Abramson --- \textit{email:} josh@stat.berkeley.edu; \textit{address:} Department of Statistics, University of California, 367 Evans Hall, Berkeley, CA 94720-3860}




\begin{abstract} 
We study the structure of the shocks for the inviscid Burgers equation
in dimension 1 when the initial velocity is given by L\'evy noise, or
equivalently when the initial potential is a two-sided L\'evy process
$\psi_0$. When $\psi_0$ is abrupt in the sense of Vigon or has bounded
variation with $\limsup_{|h| \downarrow 0} h^{-2} \psi_0(h) = \infty$, we
prove that the set of points with zero velocity is regenerative, and
that in the latter case this set is equal to the set of Lagrangian
regular points, which is non-empty. When $\psi_0$ is abrupt we show that the shock
structure is discrete. When $\psi_0$ is eroded we show that there are
no rarefaction intervals.
\end{abstract}

\maketitle

\section{Introduction}
\label{sec:intro}

Burgers introduced the equation
\[
\partial_t u + \partial_x (u^2/2) = \eps \partial_{xx}^2 u
\]
as a simple model of hydrodynamic turbulence for compressible
fluids, where the parameter $\eps > 0$ describes the viscosity
of the fluid and the solution represents the velocity of a fluid particle
located at $x$ at time $t$~\cite{burgers}.
It can be seen as a simplification of the Navier-Stokes
equation arrived at by neglecting pressure and force terms, but also
arises in other physical problems, such as the formation of the
superstructure of the universe~\cite{bert20}.

It is known that under certain conditions, as $\eps \ra 0$ the
solution converges to the unique entropy condition satisfying weak
solution of the \emph{inviscid Burgers equation}
\begin{equation}
\label{eq:inviscid_burgers}
\partial_t u + \partial_x (u^2/2) = 0 .
\end{equation}
A physical interpretation of the weak entropy condition satisfying
solution to~\eqref{eq:inviscid_burgers} is that at time zero,
infinitesimal particles are uniformly spread on the line, with initial
velocity $u( \cdot,0)$, and these particles evolve according to the
dynamics of completely inelastic shocks. That is, the velocity of a
particle changes only when the cluster of particles it is in collides
with another cluster, in which case the clusters stick together and
form a heavier cluster, with conservation mass and momentum
determining the mass and velocity of the new cluster.

There is an abundant literature on the solution to~\ref{eq:inviscid_burgers}
when the initial velocity $u(\cdot,0)$ is a random process. See for example \cite{bert1,bert2,giraud_burgers,giraud_geneology,bertoin_burgers_stable,bertoin_large_devs,burgers,bert7,bert12,bert15,bert16,bert18,bert19,bert20,l-rey,winkel_limit}.
We will investigate the solution when $u(\cdot,0)$ is a L\'evy noise,
i.e. when the \emph{potential} process $\psi_0 = (\psi_0(x))_{x \in
  \reals}$, defined by $\psi_0(x) - \psi_0(y) = \int_x^y u(z,0) \,
dz$, has stationary independent increments. In particular, we
investigate qualitative features of the shock structure of the
solution, and thus extend the work of 
Bertoin~\cite{bertoin_burgers_stable} ($\psi_0$ a stable L\'evy process with stability index $\alpha \in (1/2,2]$),
Giraud~\cite{giraud_burgers} (extensive results for the case $\alpha
\in (1/2,1)$) and Lachi\`{e}ze-Rey~\cite{l-rey} ($\psi_0$ a bounded
variation L\'evy process).

In order to explain our results, we must first discuss the general solution
to~\eqref{eq:inviscid_burgers} and some related concepts. We follow
\cite[Section $2.1$]{giraud_burgers} closely.
Suppose that $\psi_0$ has discontinuities only of the first kind and
satisfies $\psi_0(x) = o(x^2)$ as $|x| \ra \infty$. Then as $\eps \ra
0$ the unique solution of Burgers equation with viscosity $\eps > 0$
converges (except on a countable set) to a weak solution
of~\eqref{eq:inviscid_burgers}, referred to as the Hopf-Cole solution
(see~\cite{hopf,cole}). The right continuous version of this solution
is
\[
u(x,t) = t^{-1} (x-a(x,t)),
\]
where, taking the supremum over all possible arguments if necessary,
\[
a(x,t) := \arg \sup \left\{ \psi_0(y) - \mbox{$\frac{1}{2t}$}(y-x)^2 \,
  : \, y \in \reals \right\} .
\]
The function $x \mapsto a(x,t)$ is non-decreasing and right continuous
and its right continuous inverse $a \mapsto x(a,t)$ is known as the
\emph{Lagrangian function}, and gives the position at time $t$ of the
particle initially located at $a$.

A discontinuity of $x \mapsto u(x,t)$ is called a shock and occurs
when $x \mapsto a(x,t)$ jumps, i.e. when $a(x,t) \neq a(x-,t) :=
\lim_{y \uparrow x} a(y,t)$. From the point of view of the particle
description, the location of a shock corresponds to the location of a
cluster at time $t$. This cluster results from the aggregation of the
particles initially located in $[a(x-,t),a(x,t)]$; its velocity is
(according to the conservation of masses and momenta)
\[
v(x, t) = - \frac{\psi_0(a(x, t))-\psi_0(a(x-, t))}{a(x, t)-a(x-, t)} =
\mbox{$\frac{1}{2}$} \left[ u(a(x,t)) + u(a(x-,t)) \right] .
\]
The interval $[a(x-,t),a(x,t)]$ is called a \emph{shock interval} and
$x$ a \emph{Eulerian shock point}. We define the \emph{shock
  structure} of the solution at time $t$ to be the closed range of
$a(\cdot,t)$. Of particular interests are points which are not
isolated on the left or the right in that closed range, since they
represent the initial locations of particles that have not been involved in any
collisions by time $t$. We call any such point a \emph{Lagrangian
  regular point}.
Finally, we call $(x,y)$ a \emph{rarefaction} interval if $a(\cdot,t)$ stays
constant on $[x,y)$. A rarefaction interval represents an interval
where there are no fluid particles at time $t$.

Our results concern qualitative features of the shock structure, 
the regenerativity of the process $(u(x,t))_{x \in \reals}$ at points
where $u(x,t) = 0$, and the relationship between such points and the
Lagrangian regular points. For our arguments, there is no loss of
generality to assume $t=1$ -- the properties we show will be true for any $t >
0$. Thus we restrict our attention to the case $t=1$ and set $a(x) =
a(x,1)$, $u(x) = u(x,1)$ for all $x \in \reals$. 
The shock structure is then
\[
\mathcal{A} := \mathbf{cl} \{ y \in \reals \, : \, a(x) = y \text{ for some } x \in
\reals \} ,
\]
i.e.\ the closure of the range of $a(\cdot)$,
and Lagrangian regular points are the subset of points of
$\mathcal{A}$ that are neither left nor right isolated. We also define
$\mathcal{A}_0 \subset \mathcal{A}$ by
\[
\mathcal{A}_0 := \mathbf{cl} \{ x \in \reals \, : \, a(x) = x \} 
= \mathbf{cl} \{ x \in \reals \, : \, u(x) = 0 \} ,
\]
Note that both $\mathcal{A}$
and $\mathcal{A}_0$ are stationary sets when $\psi_0$ is a L\'evy
process, and since adding a drift term has no affect on the
distributions of these random sets, \emph{we will assume throughout
  that if $\psi_0$ has bounded variation then it has zero drift
  coefficient}.

To ensure that $\mathcal{A}$ is non-empty we will always assume that
$\lim_{|x| \ra \infty} x^{-2} \psi_0(x) = 0$, and in the bounded
  variation case we mostly assume that $\limsup_{|h| \downarrow 0} h^{-2}
\psi_0(h) = \infty$ to ensure that $\mathcal{A}$ has a nice
structure. Most of our results in the bounded variation case also
require a further assumption relating to overshoots at hitting times
-- see Assumption~\textbf{B} in Section~\ref{sec:hypotheses}.

In all cases we show that the Lebesgue measure of $\{ y \in \reals \,
: \, a(x) = y \text{ for some } x \in \reals \}$ is zero (see
Lemma~\ref{lem:null_set} ) and in the bounded variation case we show
that this set is closed (see Theorem~\ref{thm:Aclosed}).  For $\psi_0$
in an interesting class of unbounded variation L\'evy processes called
\emph{abrupt} L\'evy processes (see
Section~\ref{sec:abrupt_and_eroded} for a definition), we also show
that this set is closed and that moreover $\mathcal{A}$ is a discrete
set (see Theorem~\ref{thm:abrupt} and Corollary~\ref{cor:abrupt}),
extending the result of~\cite{bertoin_burgers_stable} that this is
true when $\psi_0$ is a stable process with $\alpha \in (1,2]$. A
result from~\cite{bertoin_burgers_stable} relating to Cauchy processes
is also extended to a more general class of unbounded variation
processes, the \emph{eroded} L\'evy processes (again, see
Section~\ref{sec:abrupt_and_eroded} for a definition).  For these
eroded processes, there are no rarefaction intervals.

We show that if $\psi_0$ is of unbounded variation and abrupt, or of
bounded variation and satisfying Assumption~\textbf{B}, then the
process $u = (u(x))_{x \in \reals}$ is regenerative at points $y$ such
that $u(y) = 0$, that between any two consecutive such points it must
first be positive and then negative, and that the only accumulation
points of jump times of $u$ are at such points (see
Theorem~\ref{thm:regenerative}, Theorem~\ref{thm:abrupt},
Proposition~\ref{prop:bounded} and
Theorem~\ref{thm:zeros_are_regular}). For $\psi_0$ a stable processes
with $\alpha \in (1/2,1)$, this is the main result
of~\cite{giraud_burgers}, hence our work generalizes that result to a
wider class of bounded variation processes (it is also shown
in~\cite{giraud_burgers} that for those stable processes $\mathcal{A}$
is a discrete set -- we could not generalize this result to our wider
class of bounded variation processes).  Key to proving this result is
the theory of \emph{randomized coterminal times} due to Millar (see
Section~\ref{sec:rct}), which allows us to decompose the process at $T
:= \inf \{ x \ge 0 : x \in \mathcal{A}_0 \}$, i.e. at the first
non-negative element of $\mathcal{A}_0$. The results of
Lachi\`{e}ze-Rey~\cite{l-rey} also form an indispensable part of our
arguments in the bounded variation case.

Another important result of~\cite{giraud_burgers} is that when
$\psi_0$ is a stable processes with $\alpha \in (1/2,1)$,
$\mathcal{A}_0$ is exactly equal to the set of points of $\mathcal{A}$
at which $\psi_0$ is continuous, which is in turn equal to the set of
Lagrangian regular points. We extend this result to our more general
class of bounded variation processes (see
Proposition~\ref{prop:bounded} and
Theorem~\ref{thm:zeros_are_regular}) again using the results of
Lachi\`{e}ze-Rey~\cite{l-rey}.

The rest of the paper is organized as follows. In
Section~\ref{sec:geometry} we discuss geometric interpretations of
$a(x)$ that make the proofs easier to read, and introduce the important
connection between $\mathcal{A}$ and the concave majorant of
$(\psi_0(x) - \frac{1}{2} x^2)_{x \in \reals}$. In
Section~\ref{sec:main} we present and prove all of our results, with
the exception of the proof of the regenerativity property of
$\mathcal{A}_0$ mentioned above, which we prove in
Section~\ref{sec:regenproof}.

We conclude the introduction by noting that $\mathcal{A}_0$ is the set
of fixed points of the proximal mapping for the Moreau envelope of
$\psi_0$~\cite{MR2496900,MR1491362} and thus may be of interest in
convex analysis.

\newpage

\section{Geometric Interpretations and Relation to Concave Majorants}
\label{sec:geometry}

Recall from Section~\ref{sec:intro} that
\[
a(x) = \arg \sup \left\{ \psi_0(y) - \mbox{$\frac{1}{2}$}(y-x)^2 \,
  : \, y \in \reals \right\} ,
\]
i.e. $a(x)$ is the (largest) location of the supremum of $y \mapsto
\psi_0(y) - \frac{1}{2}(y-x)^2$. One has the following geometric
interpretation: consider a realization of the initial potential
$\psi_0$ and a parabola $y \mapsto \frac{1}{2}(z-x)^2 + C$, where $C$
is chosen such that the parabola is strictly above the path of
$\psi_0$. Let $C$ decrease until this parabola touches the graph of
$\psi_0$. Then $a(x)$ is the largest abscissa of the contact points.

Now consider what happens to $a(x)$ as $x$ increases. Suppose for
example that $x < a(x)$, then the center of
the parabola will move forward, and $C$ will increase so that the
largest abscissa of the contact points between the parabola and
$\psi_0$ remains at $a(x)$. This will keep going until for some $z >
x$, the location of the largest supremum of $y \mapsto
\psi_0(y) - \frac{1}{2}(y-z)^2$ is no longer at $a(x)$, that is, the
parabola centered at $z$ passing through the point
$(a(x),\psi_0(a(x)))$ will touch $\psi_0$ again at
$(a(z),\psi_0(a(z)))$, where $a(z) > a(x)$. 
This creates a jump in $a$, and hence in $u$, at the location
$z$. The story is similar when $x > a(x)$, except that now $C$ will
decrease in order to keep the parabola touching $\psi_0$ as the center
of the parabola moves forward.

Another important geometric property of the Hopf-Cole solution relates
to concave majorants. For any $f: \reals \to \reals$, the concave
majorant of $f$ is the minimal concave function $\bar{C}_f: \reals
\to \reals \cup \{ \infty \}$ such that $\bar{C}_f(x) \ge f(x) \vee
f(x-)$ for every $x \in \reals$.

Let $\bar{C} : \reals \ra \reals$ denote the concave majorant of
$(\psi_0(x) - \frac{1}{2}x^2)_{x \in \reals}$, and denote its right
continuous derivative by $\bar{c} = \bar{C}'$. Since $\bar{c}(\cdot)$
is non-increasing, we can consider the Stieltjes measure $-d
\bar{c}$. The connection with $\mathcal{A}$ is the following.

\begin{lem}
\label{lem:contained}
For any $\psi_0$, 
\begin{enumerate}[(i)]
\item $\text{\emph{Supp}}(d \bar{c}) \subseteq \textbf{cl} \{ y \, :
  \, \exists \, x \text{ s.t. } a(x) = y \} = \mathcal{A}$;
\item $\{ y \, : \, \exists \, x \text{ s.t. } a(x) = y \} \subseteq
  \text{\emph{Supp}}(d \bar{c})$.
\end{enumerate}
Hence if $\{ y \, : \, \exists \, x \text{ s.t. } a(x) = y \}$ is
closed then $\text{Supp}(d \bar{c}) = \mathcal{A}$.
\end{lem}
\proof
(i)
Suppose first that $y \in \text{Supp}(d \bar{c})$ is isolated on both
sides in $\text{Supp}(d \bar{c})$ or is in the interior of
$\text{Supp}(d \bar{c})$. Then there exists $x \in \reals$ such that
\[
\left( \psi_0(y+z) \vee \psi_0((y+z)-) - \mbox{$\frac{1}{2}$} (y+z)^2 \right)
- \left( \psi_0(y) \vee \psi_0(y-) - \mbox{$\frac{1}{2}$} y^2 \right) 
< - x z
\]
for all $z \neq 0$. But then 
\[
\psi_0(y + z) \vee \psi_0((y+z)-) - \mbox{$\frac{1}{2}$}(y+z-x)^2
\leq \psi_0(y) \vee \psi_0(y-) - \mbox{$\frac{1}{2}$} (y-x)^2
\]
with equality only if $z = 0$. Hence
\[
a(x) = y + \arg \sup \left\{ \psi_0(y+z) - \mbox{$\frac{1}{2}$}((y+z)-x)^2 \,
  : \, z \in \reals \right\} = y + 0 = y ,
\]
and thus $y \in \mathcal{A}$.

Now suppose $y$ is not isolated in $\text{Supp}(d \bar{c})$. Then
there exists a sequence of points $\{ y_n \}_{n \ge 0}$ such that $y_n
\ra y$ with each $y_n$ either isolated on both sides in $\text{Supp}(d
\bar{c})$ or in the interior of $\text{Supp}(d \bar{c})$. Let $\{ x_n
\}_{n \ge 0}$ be such that $a(x_n) = y_n$ for each $n \ge 0$. Then
$a(x_n) \ra y$ and hence $y \in \mathcal{A}$ since $\mathcal{A}$ is
closed.  

(ii)
Suppose there exists $x$ such that $a(x) = y$. From the definition of
$a(x)$ it follows that 
\[
\psi_0(y - z) \vee \psi_0((y-z)-) - \mbox{$\frac{1}{2}$} ((y-z)-x)^2
\leq 
\psi_0(y) \vee \psi_0(y-) - \mbox{$\frac{1}{2}$} (y-x)^2
\]
for all $z \ge 0$ and
\[
\psi_0(y + z) \vee \psi_0((y+z)-) - \mbox{$\frac{1}{2}$} ((y+z)-x)^2
<
\psi_0(y) \vee \psi_0(y-) - \mbox{$\frac{1}{2}$} (y-x)^2
\]
for all $z>0$. Thus
\begin{equation}
\label{eq:lefty}
\left( \psi_0(y - z) \vee \psi_0((y-z)-) - \mbox{$\frac{1}{2}$}
  (y-z)^2 \right)
- \left( \psi_0(y) \vee \psi_0(y-) - \mbox{$\frac{1}{2}$} y^2 \right)
\leq zx
\end{equation}
for all $z \ge 0$ and
\begin{equation}
\label{eq:righty}
\left( \psi_0(y + z) \vee \psi_0((y+z)-) - \mbox{$\frac{1}{2}$}
  (y+z)^2 \right)
- \left( \psi_0(y) \vee \psi_0(y-) - \mbox{$\frac{1}{2}$} y^2 \right)
< - zx
\end{equation}
for all $z>0$. 

\eqref{eq:lefty} implies that $\bar{c}(y-) \geq -z$ and
\eqref{eq:righty} implies that $\bar{c}(y) < -z$. Hence $y \in
\text{Supp}(d \bar{c})$.
\endproof


\section{Definitions and Background Material}
\label{sec:defs}

\subsection{L\'evy processes}
\label{sec:levy}

Let $\psi_0 = (\psi_0(x))_{x \in \reals}$ be a real-valued L\'evy process.  
That is, 
$\psi_0$ has {\em c\`adl\`ag} sample paths, $\psi_0(0)=0$, 
and $\psi_0(y)-\psi_0(x)$ is independent of 
$(\psi_0(z))_{z \le x }$ with the same distribution as $\psi_0(y-x)$ for
all $x,y \in \reals$ with $x<y$. 

The L\'evy-Khintchine
formula says that for $x \geq 0$ the characteristic function of $\psi_0(x)$ is given by
$\E[e^{i \theta \psi_0(x)}] = e^{-x \Psi(\theta)}$ for $\theta \in \reals$,
where
\[
\Psi(\theta) = - i c \theta + \frac{1}{2} \sigma^2 \theta^2 +
\int_\reals(1 - e^{i \theta y} + i \theta y 1_{ \{ |y| < 1 \} } )
\, \Pi(dy) 
\]
with $c \in \reals$, $\sigma \in \reals_+$, and $\Pi$ a $\sigma$-finite
measure concentrated on $\reals \setminus \{0\}$ 
satisfying $\int_\reals (1 \wedge y^2) \, \Pi(dy) <
\infty$. We call $\sigma^2$ the {\em infinitesimal variance} of the
Brownian component of $\psi_0$ and $\Pi$ the
{\em L\'evy measure} of $X$. 

The sample paths of $\psi_0$ have bounded variation almost surely
if and only if $\sigma = 0$ and
$\int_\reals (1 \wedge |y| ) \, \Pi(dy) < \infty$. In this case $\Psi$
can be rewritten as
\[
\Psi(\theta) = - i d \theta + \int_\reals (1-e^{i \theta y} ) \, \Pi(dy).
\]
We call $d \in \reals$ the drift coefficient. Recall from the
introduction that we will assume $d=0$ throughout without affecting
our results.
For full details of these definitions see \cite{bertoin}.

\subsection{Fluctuation theory}
\label{sec:fluctuation_theory}

We will often make use of some basic results from fluctuation theory for L\'evy processes.

The first is due to Shtatland \cite{shtatland}.
If $\psi_0$ has paths of bounded variation with drift $d$, then
\begin{equation}
\label{eq:shtatland_statement}
\lim_{h \downarrow 0} h^{-1} \psi_0(h) = d \quad \text{a.s.}
\end{equation}
Since the jump times of $\psi_0$ form a countable set of
stopping times, by the strong Markov property it follows that for all $y$ such
that $\psi_0(y) \neq \psi_0(y-)$, i.e. at all jump times $y$ of $\psi_0$, we have
\begin{equation}
\label{eq:shtatland_jump_statement}
\lim_{h \downarrow 0} h^{-1} (\psi_0(y+h) - \psi_0(y)) = d \quad \text{a.s.}
\end{equation}
The counterpart of Shtatland's result when $\psi_0$
has paths of unbounded variation is Rogozin's result
\begin{equation}
\label{eq:rogozin_small_time}
\liminf_{h \downarrow 0} h^{-1} \psi_0(h) = - \infty
\quad \text{ and } \quad
\limsup_{h \downarrow 0} h^{-1} \psi_0(h) = + \infty \quad \text{a.s.}
\end{equation}
By  the strong Markov property, it again follows that for all $y$ such
that $\psi_0(y) \neq \psi_0(y-)$, i.e. at all jump times $y$ of $\psi_0$, we have
\begin{equation}
\label{eq:rogozin_jump}
\begin{split}
\liminf_{h \downarrow 0} h^{-1} (\psi_0(y+h) - \psi_0(y)) & = - \infty
\quad \text{a.s.} \quad \text{ and } \\
\limsup_{h \downarrow 0} h^{-1} (\psi_0(y+h) - \psi_0(y)) & = + \infty \quad \text{a.s.} \\
\end{split}
\end{equation}

\subsection{Hypotheses on $\psi_0$}
\label{sec:hypotheses}

We now define some hypotheses on $\psi_0$. We will always assume the
first and the second ensures that the shock structure is nice when
$\psi_0$ has paths of bounded variation. Let $\bar{C} : \reals \ra
\reals$ denote the concave majorant of $(\psi_0(x) -
\frac{1}{2}x^2)_{x \in \reals}$, and denote its right continuous derivative by
$\bar{c} = \bar{C}'$. Since $\bar{c}(\cdot)$ is non-increasing, we can
consider the Stieltjes measure $-d \bar{c}$.

\begin{hypA}
Let $\psi_0$ be such that almost surely $\lim_{|x| \ra \infty} x^{-2} \psi_0(x) = 0$.
\end{hypA}

\begin{remark}
\label{rem:hypA}
\begin{enumerate}[(i)]
\item Hypothesis~\textbf{A} implies $\bar{C}$ is finite and $\sup_{x \in
    \reals} \{ \psi_0(x) - \frac{1}{2}x^2 \} < \infty$. 
\item Hypothesis~\textbf{A}  holds for stable processes with stability
  index $\alpha \in (1/2,2]$. 
\end{enumerate}
\end{remark}

\begin{hypB}
If $\psi_0$ has paths of bounded variation then let $\psi_0$ be such that 
\begin{equation}
\label{eq:hypB}
\limsup_{h \downarrow 0} h^{-2} \psi_0(h) = + \infty
\quad \text{ and } \quad 
\liminf_{h \downarrow 0} h^{-2} \psi_0(h) = - \infty \quad \text{a.s.}
\end{equation}
\end{hypB}

\begin{remark}
\label{rem:hypB}
\begin{enumerate}[(i)]
\item If $\psi_0$ has paths of unbounded variation then
  \eqref{eq:hypB} always holds by~\eqref{eq:rogozin_small_time}.
\item If $\psi_0$ has paths of bounded variation then
  by~\eqref{eq:shtatland_statement} Hypothesis~\textbf{B} implies that $\psi_0$
  has zero drift coefficient (but we are already assuming this is true
  throughout). In fact, Bertoin et al.\
  have fully characterized which bounded variation L\'evy processes satisfy
  \eqref{eq:hypB}~\cite[Theorem $3.2$]{bertoin_powerlaw} (clearly
  it is necessary to at least have $\Pi((-\infty,0)) =  \Pi((0,\infty)) = \infty$).
\item Hypothesis~\textbf{B} holds for stable processes with stability
  index $\alpha \in (1/2,2]$.
\item Again by the strong Markov property, under Hypothesis~\textbf{B}
  it follows that for all $y$ such that $\psi_0(y) \neq \psi_0(y-)$,
\begin{equation}
\label{eq:hypB_jump}
\begin{split}
\limsup_{h \downarrow 0} h^{-2} (\psi_0(y+h) - \psi_0(y))  & = +
\infty  \quad \text{a.s.} \quad \text{ and } \\
\liminf_{h \downarrow 0} h^{-2} (\psi_0(y+h) - \psi_0(y)) & = - \infty \quad \text{a.s.} \\
\end{split}
\end{equation}
\end{enumerate}
\end{remark}



The following assumption will be necessary for the advanced results in
the bounded variation case. Recall that we have already assumed
$\psi_0$ to have zero drfit coefficient.

\begin{assB}
Suppose $\psi_0$ has paths of bounded variation. 
\begin{enumerate}[(I)]
\item Let $T = \inf_{x \ge 0} \{ \psi_0(x) - b x -
\mbox{$\frac{1}{2}$} x^2 \ge s \}$ for some $b > 0$ and $s > 0$.
Then on the set $\{ T < \infty \}$ we have $\psi_0(T) - b T - \frac{1}{2} T^2 > s$ almost surely. 
\item Let $T = \inf_{x \ge 0} \{ \psi_0(x) + b x -
\mbox{$\frac{1}{2}$} x^2 \le -s \}$ for some $b > 0$ and $s > 0$.
Then on the set $\{ T < \infty \}$ we have $\psi_0(T) + b T - \frac{1}{2} T^2 < s$ almost surely. 
\end{enumerate}
\end{assB}

\begin{remark}
\label{rem:assB}
\begin{enumerate}[(i)]
\item By quasi-continuity of L\'evy processes the conclusion still
  holds when $\psi_0(x)$ is replaced by $\psi_0(x) \vee \psi_0(x)$ in
  the definitions of $T$.
\item (II) has an equivalent time reversed version: 
let $T = \inf_{x \le 0} \{ \psi_0(x) + b x -
\mbox{$\frac{1}{2}$} x^2 \ge s \}$ for some $b > 0$ and $s > 0$.
Then $\psi_0(T) + b T - \frac{1}{2} T^2 > s$ a.s.
\end{enumerate}
\end{remark}

\subsection{Abrupt and eroded L\'evy processes}
\label{sec:abrupt_and_eroded}

A broad class of unbounded variation L\'evy processes of interest has
been defined by Vigon~\cite{abrupt}.

\begin{defn}
\label{def:abrupt}
A L\'evy process $\psi_0$ is {\em abrupt} if its paths have unbounded variation 
and almost surely for all $m$ such that $\psi_0$ has a local maximum at $m$,
\[
\liminf_{h \downarrow 0} h^{-1}(\psi_0(m-h)-\psi_0(m)) = +\infty 
\quad \text{ and } \quad
\limsup_{h \downarrow 0} h^{-1}( \psi_0(m+h) - \psi_0(m)) = -\infty.
\]
\end{defn}

\begin{remark}
\label{rem:abrupt}
A L\'evy process $\psi_0$ with paths of unbounded variation is abrupt
 if and only if 
\begin{equation}
\label{eq:abruptcondition}
 \int_0^1 x^{-1} \P\{ \psi_0(x) \in [ax,bx]\}  \, dx < \infty,
 \quad \forall a<b,
\end{equation}
(see \cite[Theorem 1.3]{abrupt}).
Examples of abrupt L\'evy processes include stable processes 
with stability parameter in the interval $(1,2]$, processes with
non-zero Brownian component, and any processes that creep upwards or
downwards. An example of an unbounded variation process that is not
abrupt is the symmetric Cauchy process, however this process will be
\emph{eroded} in the sense of Definition~\ref{def:eroded}.
\end{remark}

The following theorem describes the local behavior of an abrupt L\'evy
process at arbitrary times.  This result is an immediate corollary of
the more general result 
\cite[Theorem 2.6]{abrupt} once we use the fact that almost surely the
paths of a L\'evy processes cannot have both points of increase and
points of decrease \cite{fourati}.

\begin{thm}
\label{thm:allt}
Let $\psi_0$ be a two sided abrupt L\'evy process. Then,
almost surely for all $x \in \reals $, if
\[
\limsup_{h \downarrow 0} h^{-1} (\psi_0(x-h)-\psi_0(x-))  < \infty
 \quad \text{ and } \quad 
\limsup_{h \downarrow 0} h^{-1} (\psi_0(x+h)-\psi_0(x)) < \infty 
\]
then $\psi_0$ has a local supremum at $x$.
\end{thm}

At the other end of the scale from abrupt processes are eroded
processes, also introduced by Vigon~\cite{vigon_unpublished}.

\begin{defn}
\label{def:eroded}
A L\'evy process $\psi_0$ is {\em eroded} if its paths have unbounded variation 
and almost surely for all $m$ such that $\psi_0$ has a local maximum at $m$,
\[
\liminf_{h \downarrow 0} h^{-1}(\psi_0(m-h)-\psi_0(m)) = 0
\quad \text{ and } \quad
\limsup_{h \downarrow 0} h^{-1}( \psi_0(m+h) - \psi_0(m)) = 0.
\]
\end{defn}

Vigon~\cite[Theorem $1.4$]{vigon_unpublished} gives the following
characterization of eroded processes (the result may also be found
in~\cite[Theorem $3.11$]{lipz}). 

\begin{remark}
\label{rem:eroded}
A L\'evy process $\psi_0$ with paths of unbounded variation is eroded
 if and only if 
\begin{equation}
\label{eq:erodedcondition}
 \int_0^1 x^{-1} \P\{ \psi_0(x) \in [ax,bx]\}  \, dx = \infty,
 \quad \forall a< 0 < b,
\end{equation}
\end{remark}

\subsection{Randomized coterminal times}
\label{sec:rct}

Randomized coterminal times were introduced by Millar in order to
extend the set of times at which some sort of decomposition into two
independent processes could take place~\cite{coterminal}. Essentially
they are last exit times from randomized sets. For example, the
largest time at which the supremum of a Markov process $(\phi(x))_{x
  \ge 0}$ is achieved is the last exit time from the random
interval $[\sup_x \phi(x) , \infty)$.

In this subsection we assume $(\phi(x))_{x \ge 0}$ is a c\`adl\`ag
strong Markov process with state space $(E, \mathcal{E})$, a locally
compact metric space (in fact we will only use state space
$([0,\infty) \times \reals, \, \mathcal{B}([0,\infty) \times
\reals))$. Denote by $\mathcal{F}_x$ the sigma fields that are the
right continuous completions of the natural sigma fields
$\mathcal{F}_x^0 = \sigma \{ \phi(y), y \le x \}$, and let
$\mathcal{F} = \bigvee_{x \ge 0} \mathcal{F}_t$. Let $\theta_x$ be the
standard shift operator, so that $\phi(y)(\theta_x \omega) =
\phi(x+y)(\omega)$ for every $y \ge 0$. Recall that a \emph{random
  time} $R$ is a $[0,\infty]$-valued $\mathcal{F}$-measurable random
variable, and that a random time $T$ is a \emph{terminal time} if it
is optional and $T = x + T \circ \theta_x$ on $\{ T > x \}$. For a
random time $R$, define
\[
\begin{split}
\mathcal{F}(R+) := \big\{ F \in \mathcal{F} : & \text{ for all } x > 0,
\text{ there exists } F_x \in \mathcal{F}_x \\
& \qquad \text{ such that } F \cap
\{ R < x \} = F_x \cap \{ R < x \} \big\}. \\
\end{split}
\]

\begin{defn}
\label{def:rct}
Suppose we are given
\begin{itemize}
\renewcommand{\labelitemi}{$\cdot$}
\setlength{\itemindent}{-7pt}
\item a measure space $( A, \mathfrak{U} )$,
\item a family of terminal times $\{ T_a \}_{ a \in A }$
  such that $(a, \omega) \rightarrow T_a(\omega)$ is
  $\mathfrak{U} \times \mathcal{F}$-measurable,
\item a measurable mapping $Z$ from $(\Omega, \mathcal{F})$ to $(A, \mathfrak{U})$.
\end{itemize}
A random time $R$ is a \emph{randomized coterminal time} based on
$(A,\mathfrak{U})$, $\{ T_a \}_{a \in A}$, $Z$ if
\begin{enumerate}[(I)]
\item for each $x \ge 0$ there is an $\mathcal{F}_x$-measurable
  $A$-valued random variable $Z_x$ such that $Z = Z_x$ on the set $\{
  R \le x \}$,
\item for each $0 \le y < x$ there exists $B(y,x) \in \mathcal{F}_t$
  such that 
\[
\{ y \le R < x \} = B(y,x) \cap \{ T_{Z(\omega)} (\theta_x \omega) =
+\infty \} .
\]
\end{enumerate}
Note that by (I) the $Z$ in (II) can be replaced by $Z_x$.
\end{defn}

\begin{example1}
Suppose $\lim_{x \ra \infty} \psi_0(x) = -\infty$, then
$R = \arg \sup \{ \psi_0(x) \, : \, x \ge 0 \}$ is a randomized
coterminal time with $(A, \mathfrak{U}) = (\reals,
\mathcal{B}(\reals))$, $T_a = \inf \{ x > 0 \, : \, \psi_0(x) \vee
\psi_0(x-) \geq a \}$, $Z = \sup_{x \ge 0} \psi_0(x)$ and $Z_x =
\sup_{y \leq x} \psi_0(y)$. Property (I) is immediate since if the
supremum occurs before $x$, then it is equal to the supremum
attained by $\psi(y)$ on $[0,x]$, and to see that property (II) holds
note that
\[
\begin{split}
\{ y < R \leq x \} 
& = \{ \text{ the supremum of $\psi_0$ occurs in $(y,x]$ } \} \\
& = \{ \text{ $\psi_0$ goes at least as high in $(y,x]$  as it did
  before time $y$, } \color{white} \} \\
& \quad \{ \color{black} \text{ and never after $x$ goes as high as it
did during $(y,x]$ } \} \\
& = \{ Z_{yx} \leq Z_x \} \, \{ T_{Z_x(\omega)}( \theta_x \omega ) = +
\infty \} \, ,
\end{split}
\]
where $Z_{yx} = \sup_{y < w \leq x} \psi_0(w)$, so $B(y,x) = \{ Z_{yx}
\leq Z_x \}$ here.
\end{example1}

The following result is~\cite[Theorem $3.4$]{coterminal} and
essentially says that for a randomized coterminal
time $R$, conditional on $Z=z$, the post $R$ process is `just' the
original process conditioned on $\{ T_z = + \infty \}$, and is still
Markovian. Note that $Z$ is $\mathcal{F}(R+)$ measurable by (I).

\begin{thm}
\label{thm:rct_markov}
Let $(\phi(x))_{x \ge 0}$ be a Hunt process, and $R$ a randomized
coterminal time based on $(A,\mathfrak{U})$, $\{ T_a \}_{a \in A}$, $Z$. Then for
bounded Borel $f$,
\[
\E \left( f( \phi(R+x) ) | \mathcal{F}((R+y)+) \right)
=
\mbox{$\int$}  f(b) \, H_{x-y}(Z ; \phi(R+y), db), \quad 0 < y < x,
\]
where $H_x(z; a, db) := \P^a(\phi(x) \in db) \P^b(T_z = \infty) /
\P^a(T_z = \infty)$.
\end{thm}

The next result is a combination of~\cite[Proposition
$5.4$]{coterminal} and~\cite[(a) following Proposition
$5.4$]{coterminal}, where we have trivially extended the state
space to include a deterministic element as well as a L\'evy
process. It gives conditions under which, conditionally given $Z$ and
$\phi(R)$, the post $R$ process is independent of
$\mathcal{F}(R+)$. The proof relies on the zero-one property of L\'evy
processes at local maxima or jump times.

\begin{prop}
\label{prop:local_max_and_jump_rct}
Let $(\psi(x))_{x \ge 0}$ be a L\'evy process and let $\phi(x) =
(x,\psi(x))$ for $x \geq 0$. Let $R$ be a randomized
coterminal time for $(\phi(x))_{x \ge 0}$ based on $(A,\mathfrak{U})$, $\{ T_a \}_{a \in A}$, $Z$. Suppose that
$\P( R \text{ is the time of a local maximum of } \psi) = \P( R < \infty
)$ or $\P( R \text{ is a jump time of } \psi) = \P( R < \infty)$. Then
conditional on $Z$ and $\phi(R)$, the post $R$ process is independent
of $\mathcal{F}(R+)$, and it is Markov with transitions $H_x(Z ; a, db )$.
\end{prop}

\section{Main Results}
\label{sec:main}

In this section we first present results in a general setting and then
treat processes with paths of unbounded and bounded variation
separately. Recall from Remark~\ref{rem:hypB} that
Hypothesis~\textbf{B}  is automatically satisfied when $\psi_0$ has
paths of unbounded variation, and when $\psi_0$ has
paths of bounded variation then by assumption the
drift coefficient of $\psi_0$ is zero.

\subsection{Unbounded and Bounded Variation}
\label{sec:universal}

\begin{lem}
\label{lem:null_set}
Let $\psi_0$ be a two-sided L\'evy process satisfying 
satisfying Hypotheses~\textbf{A} and \textbf{B}.
The Lebesgue measure of $ \{ y \in \reals \, : \, a(x) = y \text{ for
  some } x \in \reals \}$ is zero a.s.
\end{lem}
\proof
  By application of Fubini's theorem and stationarity it suffices to
  show that $\P ( a(x) = 0 \text{ for some } x \in \reals ) = 0$.
  Suppose there exists $x > 0$ such that $a(x) = 0$, then $\limsup_{h
    \downarrow 0} h^{-1} \psi_0(h) \leq - x < 0$. But this happens
  only on an event of probability zero by
  \eqref{eq:shtatland_statement} or \eqref{eq:rogozin_small_time} for
  processes with bounded or unbounded variation respectively.

Similarly, if there exists $x < 0$
such that $a(x) = 0$, then $- \liminf_{h \uparrow 0} h^{-1} \psi_0(h)
\leq - x < 0$. But this happens only on an event of probability zero
by the time reversed versions of
\eqref{eq:shtatland_statement} and \eqref{eq:rogozin_small_time}.

Finally, if $a(0) = 0$ then $\limsup_{h \downarrow 0} h^{-2} \psi_0(h)
\leq 1 < \infty$, which by~\eqref{eq:hypB} only occurs on a set of
probability zero.
\endproof

\begin{lem}
\label{lem:positive_then_negative}
Let $\psi_0$ be a two-sided L\'evy process satisfying Hypotheses
\textbf{A} and \textbf{B}. Let
$a^-,a^+ \in \mathcal{A}_0$, i.e. let $u(a^-) = u(a^+) = 0$, with $a^- <
a^+$ and $u(x) \neq 0$ for $a^- < x < a^+$. Then there exists $a^- <
x_0 < a^+$ such that $u(x) > 0$ for all $a^- < x < x_0$ and $u(x) < 0$
for all $x_0 < x < a^+$. 
\end{lem}
\proof
From \eqref{eq:hypB} we have $a(0) > 0$ a.s. and hence
\[
\begin{split}
a_0^+ & := \inf \{ x \ge 0 : x \in \mathcal{A}_0 \} > 0 \quad
\text{a.s., and} \\
a_0^- & := \inf \{ x \ge 0 : - x \in \mathcal{A}_0 \} > 0 \quad
\text{a.s.} , \\
\end{split}
\]
where we have applied time reversal to get the second inequality.

By stationarity, it suffices to show that the claim is true for
$a^+=a_0^+$ and $a^-=a_0^-$, since any almost sure behaviour of $u$
over the interval $(a_0^-,a_0^+)$ must be shared by $u$ over
$(a^-,a^+)$ for any two consecutive members $a^- < a^+$ of
$\mathcal{A}_0$.
Define $$R := \inf \left\{ y \geq 0 \, : \, \psi_0(y-x) \vee \psi_0((y-x)-) - \psi_0(y) \vee \psi_0(y-) \leq \mbox{$\frac{1}{2}$} x^2
\text{ for all } x > 0 \right\}.$$
Since $R$ is a stopping time, \eqref{eq:hypB} implies that $R < a(R)$ a.s.,
and hence $u(R) < 0$ a.s. So we cannot have
$u(x) \ge 0$ for all $x \in (a_0^-,a_0^+)$. A time reversal argument
then implies that we cannot have $u(x) \le 0$ for all $x \in
(a_0^-,a_0^+)$.

Since $a(x)$ is non-decreasing, $u$ has only downwards jumps, and
thus $u$ cannot go from being negative to positive without passing
through zero. Hence we have the a.s. existence of the $x_0$ in the claim.
\endproof

The proof of the following Theorem is in Section~\ref{sec:regenproof}.

\begin{thm}
\label{thm:regenerative}
Let $\psi_0$ be an abrupt two-sided L\'evy process with paths of
unbounded variation satisfying Hypothesis \textbf{A} or a two-sided
L\'evy process with paths of bounded variation satisfying Hypotheses~\textbf{A}
and \textbf{B} and Assumption~\textbf{B}. Define
\[
T := \inf \{ x \geq 0 \, : \, x \in \mathcal{A}_0 \} \, .
\]
Then $(\psi_0(T+x)-\psi_0(T))_{x \ge 0}$ is independent of
$(\psi_0(T-x))_{x \ge 0}$.  As a consequence, the processes $(u(T +
x))_{ x \ge 0 }$ and $(u(T - x))_{x \ge 0}$ are independent and
$\mathcal{A}_0$ is a regenerative set.
\end{thm}

It is important to relate $\mathcal{A}_0$ to the set of Lagrangian
regular points, when such points exist. As we shall see in
Theorem~\ref{thm:zeros_are_regular}, when $\psi_0$ is a two-sided
L\'evy processwith paths of bounded variation satisfying
Hypotheses~\textbf{A} and \textbf{B} and Assumption~\textbf{B},
$\mathcal{A}_0$ is exactly equal to the set of Lagrangian regular
points.

\subsection{Unbounded variation}
\label{sec:unbounded}

\begin{lem}
\label{lem:unbounded_continuous}
Let $\psi_0$ be a two-sided L\'evy process satisfying Hypothesis
\textbf{A} with paths of unbounded variation. Then $\psi_0$ is
continuous at every point in the set $\{ y \in \reals \, : \, a(x) = y
\text{ for some } x \in \reals \}$ a.s.
\end{lem}

\begin{proof} 
From \eqref{eq:rogozin_jump} and a time reversal argument,
it follows that almost surely for every $y$ such that $y$ is a jump
time of $\psi_0$, i.e. $\psi_0(y) \neq \psi_0(y-)$,
\[
\limsup_{h \downarrow 0} h^{-1} (\psi_0(y+h) - \psi_0(y)) = + \infty
\quad \text{ and } \quad
\limsup_{h \downarrow 0} h^{-1} (\psi_0(y-h) - \psi_0(y-)) = + \infty .
\]
If $y = a(x)$ or $y = a(x-)$ for some $x$, and if $y$ is such that say
$\psi_0(y) > \psi_0(y-)$, then for every $h > 0$ we have
\[
\psi_0(y) - \mbox{$\frac{1}{2}$} (x-y)^2 \geq \psi_0(y + h) -
\mbox{$\frac{1}{2}$}  (x-y-h)^2 .
\]
Therefore we would have
\[
\limsup_{h \downarrow 0} h^{-1} (\psi_0(y+h) - \psi_0(y) ) \leq y-x < \infty ,
\]
which is impossible, except on an event with probability zero. The
case of a negative jump is similar, working now at the left of the
jump.
\end{proof}


\begin{cor}
\label{cor:local_max}
Let $\psi_0$ be a two-sided abrupt L\'evy process satisfying
Hypothesis \textbf{A}. Then $\psi_0$ has a local
supremum at every point in $\{ y \in \reals \, : \, a(x) = y
\text{ for some } x \in \reals \}$ a.s.
\end{cor}

\begin{proof}
Take any point $y \in \mathcal{A}$ and let $x$ be such that $a(x) =
y$. Then for every $z \ge 0$ we have
\[
\psi_0(y+z) \vee \psi_0((y+z)-) 
-
\psi_0(y) \vee \psi_0(y-) 
\leq
\mbox{$\frac{1}{2}$} ( z - (x-y))^2
-
\mbox{$\frac{1}{2}$} (x-y)^2 . 
\]
Recall from Lemma~\ref{lem:unbounded_continuous} that $\psi_0$ is
a.s. continuous at $y$, thus almost surely
\[
\limsup_{h \downarrow 0} h^{-1} (\psi_0(x-h)-\psi_0(x-))  \leq (x-y)
 \]
and
\[
\limsup_{h \downarrow 0} h^{-1} (\psi_0(x+h)-\psi_0(x))  \leq - (x-y)
\, .
\]
Theorem~\ref{thm:allt} then implies that $\psi_0$ must have a local
supremum at $y$.
\end{proof}

For abrupt L\'evy processes, the shock structure is discrete.

\begin{thm}
\label{thm:abrupt}
Let $\psi_0$ be a two-sided abrupt L\'evy process satisfying
Hypothesis \textbf{A}. Then $\mathcal{A}$ is a discrete set a.s.
\end{thm}

\begin{proof}
Because the random set $\mathcal{A}$ is stationary (i.e. its law is invariant by translation), we
have to prove that $ \# \{ [1,2] \cap \mathcal{A} \} < \infty$ a.s. It is easy to verify that the probability
that $a(x) \in [1,2]$ for some $x$ with $|x| > n$ goes to zero as $n
\ra \infty$, so it suffices in fact to establish that for each fixed
$n$ larger than some $n_0$,
\begin{equation}
\label{eq:finite_card}
\# \{ a(x) \in [1,2] \, : \, |x| \leq n \} < \infty \quad \text{a.s.}
\end{equation}

Suppose first that $\E | \psi_0(1) | < \infty$. 
Let $n_0$ be large enough such that $| \E \psi_0(1) | < 2n_0$.  
Now, if a point $y \in [1,2]$ can be expressed as $y = a(x)$ for some
$x \in [-n,n]$, then
\[
\psi_0(y \pm h) \vee \psi_0((y \pm h)-)< \psi_0(y) \vee \psi_0(y-) + 2n h \quad \text{ for every } h \in (0,2] ,
\]
and since by Lemma~\ref{lem:unbounded_continuous} $\psi_0$
is continuous at $y$,
\[
\psi_0(y \pm h) < \psi_0(y) + 2n h \quad \text{ for every } h \in (0,2] .
\]

Defining the Lipschitz majorant of $\psi_0$
equivalently to how the Lipschitz minorant of $\psi_0$ is defined in
\cite{lipz}, we have that if $\psi_0$ is continuous at $y$, then $y$ is in
the contact set of the $2n$-Lipschitz majorant of $\psi_0$ if and only if
$\psi_0(y \pm h) - \psi_0(y) < 2nh$ for all $h \in \reals$ 
(note that the existence of the Lipschitz majorant follows from our
assumption that $| \E \psi_0(1) | < 2n_0 \leq 2n$).  
From \cite[Theorem 3.8]{lipz} it
follows that there are only finitely many such contact points $y$ in
the interval $[1,2]$ almost surely. Suppose it were the case that
\[
\P( \# \{ a(x) \in [1,2] \, : \, |x| \leq n \} = \infty ) > 0.
\]
Then with positive probability there would exist $y_\ell, y_r \in
[1,2]$ such that $y_\ell < y_r$ and $\# \{ a(x) \in [y_\ell, y_r ] \,
: \, |x| \leq n \} = \infty $. Moreover, by the law of large numbers
applied to the left and to the right, with positive probability there
would exist such a pair with both $y_\ell$ and $y_r$ in the contact
set of the $2n$-Lipschitz majorant of $\psi_0$. If both $y_\ell$ and $y_r$
were in the contact set of the majorant, then every element of the
infinite set $\{ a(x) \in [y_\ell, y_r ] \, : \, |x| \leq n \}$ would
also be in the contact set of the majorant, but that is an event with
zero probability. 
Hence $ \# \{ a(x) \in [1,2] \, : \, |x| \leq n \} < \infty $ a.s.

Now remove the assumption that $\E | \psi_0(1) | < \infty$.
For each $N \in \nats$ define the two-sided L\'evy process $\tilde{\psi_0}^N$ by
\[
\tilde{\psi_0}^N(x) = 
\left\{
\begin{split}
\psi_0(x) - \sum_{ \substack{ 0 \leq y \leq x : \\ \psi_0(y) \neq \psi_0(y-)} } (\psi_0(y) - \psi_0(y-)) 1_{|\psi_0(y) -
  \psi_0(y-)| > N} \quad & \text{ for } x \geq 0 \\
\psi_0(x) + \sum_{ \substack{ 0 \leq y \leq x : \\ \psi_0(y) \neq \psi_0(y-)} } (\psi_0(y) - \psi_0(y-)) 1_{|\psi_0(y) -
  \psi_0(y-)| > N} \quad & \text{ for } x < 0 \\
\end{split}
\right\}
\]
so that $\tilde{\psi_0}^N$ is identical to $\psi_0$ but with all the
jumps of magnitude greater than $N$ removed. Let $\tilde{\mathcal{A}}^N$ be
defined in the same way that $\mathcal{A}$ is for the original process
$\psi_0$. Since $\E | \tilde{\psi_0}^N | < \infty$ the above arguments
imply that $\tilde{\mathcal{A}}^N \cap [1,2]$ is a finite set almost surely for
every $N$.

From the fact that $\Pi(N,\infty) < \infty$ for every $N \in \nats$,
and the hypothesis that $\psi_0(x) = o(x^2)$ a.s. as $|x| \ra \infty$, it
follows that almost surely there exists a random $\tilde{N} \in \nats$
such that $\mathcal{A} \cap [1,2] = \tilde{\mathcal{A}}^{\tilde{N}} \cap [1,2]$.
Hence $\mathcal{A} \cap [1,2]$ is a finite set almost surely. 
\end{proof}

\begin{cor}
\label{cor:abrupt}
Let $\psi_0$ be a two-sided abrupt L\'evy process satisfying
Hypothesis \textbf{A}. Then $\{ y \in \reals \, : \, a(x) = y \text{
  for some } x \in \reals \}$ is closed a.s.
\end{cor}

The proof of the following theorem closely follows the proof of
\cite[Theorem 5]{bertoin_burgers_stable} with Cauchy processes replaced by
eroded processes.

\begin{thm}
\label{thm:eroded}
Let $\psi_0$ be a two-sided eroded L\'evy process satisfying
Hypothesis \textbf{A}.  Then with probability one there
are no rarefaction intervals.
\end{thm}
\proof
Recall from Lemma~\ref{lem:unbounded_continuous} that jump times of
$\psi_0$ do not belong to $\mathcal{A}$ almost surely. Now suppose $(x,x')$ is a
rarefaction interval, that is $a( \cdot)$ stays constant on $[x,x')$;
denote its value by $y$. As $y$ is not a jump time of $\psi_0$, we have for
all $h>0$,
\[
\begin{split}
\psi_0(y) - \mbox{$\frac{1}{2}$} (x-y)^2 & \geq \psi_0(y-h) -
\mbox{$\frac{1}{2}$}  (x-y+h)^2 , \\
\psi_0(y) - \mbox{$\frac{1}{2}$} (x'-y)^2 & \geq \psi_0(y+h) - \mbox{$\frac{1}{2}$} (x'-y-h)^2 . \\
\end{split}
\]
We deduce that 
\[
\begin{split}
\liminf_{h \downarrow 0} h^{-1} (\psi_0(y) - \psi_0(y-h)) & \geq y - x
, \\
 \limsup_{h \downarrow 0} h^{-1} (\psi_0(y+h) - \psi_0(y)) & \leq y - x' . \\
\end{split}
\]
Since $x < x'$, we can find a rational number $q \in (y-x', y -
x)$. Then $y$ is the location of a local maximum of $(\psi_0^{(q)}(x))_{x
  \in \reals}$, where $\psi_0^{(q)}(x) := \psi_0(x) - qx$, and moreover
\begin{equation}
\label{eq:eroded}
\liminf_{h \downarrow 0} h^{-1} (\psi_0^{(q)}(y) - \psi_0^{(q)}(y+h)) > 0 .
\end{equation}
On the other hand, the family $\left( \psi_0^{(s)}, s \in \mathbb{Q}
\right)$ is a countable family of eroded processes. For each of these
processes, with
probability one, for any $s \in \mathbb{Q}$ and any location $\mu$ of
a local maximum for $\psi_0^{(s)}$,
\[
\liminf_{h \downarrow 0} h^{-1} (\psi_0^{(q)}(\mu) - \psi_0^{(q)}(\mu+h))  = 0.
\]
We conclude that~\eqref{eq:eroded} is impossible, except on an
event of probability zero, and therefore almost surely there are no
rarefaction intervals.
\endproof

\subsection{Bounded variation}
\label{sec:bounded}

\begin{thm}
\label{thm:regulars_exist}
Let $\psi_0$ be a two-sided L\'evy process satisfying
Hypothesis \textbf{A} with paths of bounded
variation.
Suppose zero is regular for $[0,\infty)$ and $(-\infty,0]$ for
$(\psi_0(x))_{x \ge 0}$, then a.s. Lagrangian regular points exist.
\end{thm}
\proof
We shall prove that the time of the maximum of $(\psi_0(x))_{0
  \leq x \leq 1}$ has positive probability of being Lagrangian
regular, and as pointed out by Bertoin in a comment before the proof
of Theorem 3 of~\cite{bertoin_burgers_stable}, it is easy to deduce
from this fact that Lagrangian regular points exist with probability
one. This is because of stationarity and the asymptotic independence
of the events $A_0$ and $A_n$ as $n \ra \infty$, where $A_n :=\{ \arg
\sup_{n \le x \le n+1} \psi_0(x) \text{ is Lagrangian regular} \}$ for
$n \ge 0$.

Let $\mu$ be the almost surely unique location of the maximum of 
$(\psi_0(x))_{0 \leq x \leq 1}$. It follows from the concave majorant theory
of Pitman and Uribe-Bravo~\cite{pitmanbravo} that $\mu \in (0,1)$, and
that if $\bar{B} :[0,1] \ra \reals$ denotes the concave majorant of
$(\psi_0(x))_{0 \leq x \leq 1}$ then its derivative $\bar{b} = \bar{B}'$ is
continuous at $\mu$ and 
\begin{equation}
\label{eq:bar_b}
 \bar{b} (\mu + h) < \bar{b} ( \mu ) = 0 < \bar{b} (\mu-h)
\end{equation}
for every sufficiently small $h > 0$.

The rest of the argument is exactly as in the proof of Theorem 3
of~\cite{bertoin_burgers_stable}.

\eqref{eq:bar_b} implies that the support of the Stieltjes measure 
$-d \bar{b}$ contains $\mu$, and more precisely $\mu$ is neither
isolated to the left nor to the right in $\text{Supp}(d
\bar{b})$. Pick any $y \in \text{Supp}(d \bar{b})$ arbitrarily
close to $\mu$. Clearly, the graph of $\bar{B}$ touches that of 
$(\psi_0(x))_{0 \leq x \leq 1}$ at $y$, so we must have $\bar{B} (y) = \psi_0(y)$ or
$\bar{B} (y) = \psi_0(y-)$. In both cases, $y$ is the location of a maximum
of $x \ra \psi_0(x) - \bar{b} (y) x$ on $[0,1]$, and a fortiori $y$ is
then the unique location of the maximum of 
$x \ra \psi_0(x) - \frac{1}{2} (y- \bar{b}(y) -x)^2$ on $[0,1]$.
Plainly, $\mu$ is also the unique location of the maximum of 
$x \ra \psi_0(x) - \frac{1}{2} (\mu -x)^2$ on $[0,1]$.
Because $\psi_0(\mu) > \max (\psi_0(0),\psi_0(1))$, there is a positive probability that
the preceding two maxima are global (i.e.\ on $\reals$) and not only
local (i.e.\ on $[0,1]$). We conclude that with positive probability,
$\mu \in \mathcal{A}$ and is neither isolated on its right nor on its left, and
therefore is a Lagrangian regular point.
\endproof

The next two results are due to Lachi\'eze-Rey~\cite[Theorem
$4.3$, Proposition $5.3$]{l-rey} and 
allow us to find the behaviour of
$\psi_0$ around points of $\mathcal{A}$ in Proposition~\ref{prop:bounded}.

\begin{thm}
\label{thm:l-rey}
Let $\psi_0$ be a two-sided L\'evy process with paths of bounded
variation.
Let $\bar{C} :[0,1] \ra \reals$ denote the concave majorant of
$(\psi_0(x) - \frac{1}{2}x^2)_{x \in \reals}$ and denote its derivative by $\bar{c} =
\bar{C}'$.
Then for all $a \in \mathcal{A}$, $a$ is left isolated (resp.\ right
isolated) in $\mathcal{A}$ if $\bar{c}(a-) \neq -a$
(resp. $\bar{c}(a) \neq -a$).
\end{thm}

\begin{prop}
\label{prop:l-rey}
Let $\psi_0$ be a two-sided L\'evy process with paths of bounded
variation. Suppose $y \in \mathcal{A}$ and
$x$ is such that $a(x) = y$. Then almost surely if
$x < y$ then $\psi_0(y-) < \psi_0(y)$ and if
$x > y$ then $\psi_0(y-) > \psi_0(y)$.
\end{prop}

\begin{prop}
\label{prop:bounded}
Let $\psi_0$ be a two-sided L\'evy process satisfying
Hypothesis \textbf{A} with paths of bounded
variation.
Then for every $y \in \mathcal{A}$, almost surely
\begin{enumerate}[(i)]
\item if $\psi_0(y-) < \psi_0(y)$ then
  $y$ is left isolated in $\mathcal{A}$, and \newline
         if $\psi_0(y-) > \psi_0(y)$ then
  $y$ is right isolated in $\mathcal{A}$;
\item if $y$ is Lagrangian regular then $\psi_0$ is continuous at $y$;
\item if $y \neq a(y)$ then $y$ is isolated in $\mathcal{A}$.
\end{enumerate}
\end{prop}
\proof
(i) 
Suppose $\psi_0(y-) < \psi_0(y)$. Then $y$ will be left isolated
in the support of the Stieltjes measure $- d \bar{c}$, and hence will
be left isolated in $\mathcal{A}$ by Lemma~\ref{lem:contained}. The
argument is similar for the case $\psi_0(y-) > \psi_0(y)$.

(ii) 
Suppose $\psi_0$ is not continuous at $y$. Then either $\psi_0(y-) <
\psi_0(y)$ or $\psi_0(y-) > \psi_0(y)$ and hence (i) implies that
$y$ cannot be Lagrangian regular. 

(iii)
By hypothesis, for any $x$ such that $a(x)
= y$, we must have $x > y$ or $x < y$. Suppose $x < y$. 
Proposition~\ref{prop:l-rey} implies that $\psi_0(y-) <
\psi_0(y)$ and thus $y$ will be left isolated in
$\mathcal{A}$ by (i).  
Moreover, $a(x) = y$ implies that $- \bar{c}(y) \leq x < y$, thus $y$
will be right isolated in $\mathcal{A}$ by Theorem~\ref{thm:l-rey}.
The argument is similar in the alternative case $x > y$.
\endproof 

Proposition~\ref{prop:bounded} (iii) immediately leads to the
following corollary.

\begin{cor}
\label{cor:bounded}
Let $\psi_0$ be a two-sided L\'evy process satisfying
Hypothesis \textbf{A} with paths of bounded
variation.
Then the set $\{ x \in \reals \, : \, a(x) = x \}$ is closed a.s.
\end{cor}

Theorem~\ref{thm:l-rey} also allows us to prove the following two theorems.

\begin{thm}
\label{thm:Aclosed}
Let $\psi_0$ be a two-sided L\'evy process with paths of bounded
variation satisfying Hypothesis~\textbf{A}, Hypothesis~\textbf{B} and
Assumption~\textbf{B}(I). 
Then the set $\{ y \in \reals \, : \, a(x) = y \text{ for some } x \in
\reals \}$ is closed a.s. and hence is equal to $\mathcal{A}$ a.s.
\end{thm}
\proof
Since in the definition of $a(x)$ we take the supremum over all possible
$\arg \sup$s, we have that
\begin{equation}
\label{eq:Aclosed}
\bar{c}(y-) > \bar{c}(y+h) \, \, \forall \, \, h > 0 
\quad \Longleftrightarrow \quad 
\, \exists \, \, x \text{ s.t. } a(x) = y \, \, .
\end{equation}
Suppose $y$ is a right accumulation point of the set $\{ y \in \reals
\, : \, a(x) = y \text{ for some } x \in \reals \}$ so that there
exists a sequence $\{ y_n \}_{n \in \nats}$ with $y_n \downarrow y$
and $\bar{c}(y_n-) > \bar{c}(y_n+h)$ for all $h>0$ and hence $y \in \{
y \in \reals \, : \, a(x) = y \text{ for some } x \in \reals \}$.

Now suppose $y$ is a left but not a right accumulation point. Then by
Lemma~\ref{lem:contained} there exists $\hat{y}$ such that $\bar{c}(y+h) =
\bar{c}(y+)$ for all $0 \leq h < \hat{y} - y$ and such that
\begin{equation}
\label{eq:yhat}
\psi_0(\hat{y}) \vee \psi_0(\hat{y}-) - \frac{1}{2} \hat{y}^2 = \bar{C}(\hat{y})
\end{equation}
i.e.\ $\hat{y}$ is the next contact point after $y$ for the concave
majorant of $(\psi(x) - \frac{1}{2} x^2)_{x \in \reals}$). 

Take any $q \in \mathbb{Q}$ such that $y < q < \hat{y}$.
Let $\bar{C}^q : (-\infty, q] \ra \reals$ be
the concave majorant of $\left( \psi(x) - \frac{1}{2} x^2 \right)_{ x
  \leq q}$ and let $\bar{c}^q$ be its right continuous derivative,
which will agree with $\bar{c}$ on the set $(-\infty,y)$. Define
\[
E^q : = \{ x \leq q \, : \, \bar{c}^q(x-) = -x \} \, . 
\]
Since $x \in E^q$ implies that at least one of $\limsup_{h \downarrow 0} h^{-1} (
\psi_0(x+h) - \psi_0(x) )$ or $\limsup_{h \downarrow 0} h^{-1} (
\psi_0(x-h) - \psi_0(x-) )$ is finite, \eqref{eq:hypB} and Fubini
imply that $E^q$ has measure zero almost surely. Also, by
Theorem~\ref{thm:l-rey} we know that $y \in E^q$ a.s.\ since $y$ is not
isolated on the left.

The outline of the rest of the argument is as follows. For each $x \in
E^q$ we will define a random time that essentially is the first time
the process $ ( \psi_0(q+z) - \mbox{$\frac{1}{2}$} (q+z)^2 )_{z \ge
  0}$ is greater than or equal to the line extending out from $x$ with
slope $-x$, i.e. the same slope as the concave majorant at $x$. This
time should be the next time the process $(\psi_0(y) - \frac{1}{2} y^2
)_{y \in \reals} $ meets its concave majorant after $y$, but using
Assumption~\textbf{B} we show that it goes strictly above that line
at that time, which leads to a contradiction.

For every $x \in E^q$ define
\[
T^q(x) := \inf \left\{ z \ge 0 \, : \, \psi_0(q+z) \vee
\psi_0((q+z)-) - \mbox{$\frac{1}{2}$} (q+z)^2 \geq \psi_0(x) \vee
\psi_0(x-) - x \left( (q-x) + z \right) \right\} \, ,
\]
and note that almost surely $T^q(x) > 0$ for every $x \in
E^q$ since by \eqref{eq:hypB} $\bar{C}(q) > \psi_0(q) \vee
\psi_0(q-) - \frac{1}{2} q^2$ a.s. Also, by definition $T^q(y) = \hat{y}$.

By Assumption~\textbf{B}(I) and the fact that $E^q$ has measure zero
a.s. it follows that a.s.
\begin{equation}
\label{eq:yhat2}
\psi_0(T^q(x)) - \mbox{$\frac{1}{2}$} (q+T^q(x))^2 \geq \psi_0(x) \vee \psi_0(x-) - y \left( (q-x) + T^q(x) \right)
\end{equation}
for every $x \in E^q$ such that $T^q(x) < \infty$. But $y \in
E^q$ and $T^q(y) = \hat{y}$ a.s. hence \eqref{eq:yhat} and
\eqref{eq:yhat2} would imply that $\hat{y} = \infty$, and thus $y$ cannot be
as assumed a left accumulation point and isolated on the right. Since
we have shown the points not isolated on the right are included in the
set $y \in \{ y \in \reals \, : \, a(x) = y \text{ for some } x \in
\reals \}$, this concludes the proof.
\endproof

\begin{thm}
\label{thm:zeros_are_regular}
Let $\psi_0$ be a two-sided L\'evy process with paths of bounded
variation satisfying Hypothesis~\textbf{A}, Hypothesis~\textbf{B} and
Assumption~\textbf{B}(II). 
Then for every $y \in \reals$, $y = a(y)$ if and only if $y$ is a Lagrangian
regular point. Hence $\mathcal{A}_0$ is exactly the set of Lagrangian
regular points.
\end{thm}
\proof
Suppose $y = a(x)$  is a Lagrangian regular point, which from
Theorem~\ref{thm:l-rey}, is possible only if $\bar{c}(y-) = -y =
\bar{c}(y)$. Since $y$ is isolated neither on the left or the right in
$\mathcal{A}$, Lemma~\ref{lem:contained} implies that
\[
 \bar{c} (y + h) < \bar{c} ( y ) = -y < \bar{c} (y-h)
\]
for every $h > 0$. Thus
\[
\left( \psi_0(y) \vee \psi_0(y-) - \mbox{$\frac{1}{2}$} y^2 \right) - y s 
< 
\left( \psi_0(y+s) \vee \psi_0((y+s)-) - \mbox{$\frac{1}{2}$} (y+s)^2 \right)
\]
for all $s \neq 0$. Rearranging, we see that
\[
 \psi_0(y+s) \vee \psi_0((y+s)-) - \psi_0(y) \vee \psi_0(y-) - \mbox{$\frac{1}{2}$} s^2 > 0
\]
for all $s \neq 0$. It follows that $y = a(y)$.

 Conversely, suppose that $y = a(y)$.
If $y$ is right isolated in $\mathcal{A}$, then there exists $\hat{y} > y$
with $a(\hat{y}) = y$, and hence by Proposition~\ref{prop:l-rey} $y$ is the
time of a negative jump of $\psi_0$. However, this would imply that $y
\notin \mathcal{A}$ by the time reversed version of~\eqref{eq:hypB_jump}, and hence $y$ is not
right isolated in $\mathcal{A}$ a.s.

If $y$ is left isolated in $\mathcal{A}$ we do not yet know that there
necessarily exists an $\hat{y}$ such that $\hat{y} < y$ and $a(\hat{y}) = y$, because
although $y=a(y)$ implies that $\arg \sup \{ \psi_0(x) - \frac{1}{2}
(x-y)^2 : x \in \reals\} = y$, the supremum may not be achieved at a unique
point. Once we have shown that the supremum is unique a.s.\
a similar argument to the right isolated case above would show that
$y$ is not left isolated in $y$ a.s. and hence that $y = a(y)$ implies
that $y$ is a Lagrangian regular point.

Suppose that $\arg \sup \{ \psi_0(x) - \frac{1}{2} (x-y)^2 : x \in \reals \} = y$ and
there exists $\hat{y} < y$ such that $\psi_0(\hat{y}) \vee \psi_0(\hat{y}-) -
\frac{1}{2} (\hat{y}-y)^2 = \psi_0(y) \vee \psi_0(y-)$, i.e.\ suppose that
the supremum is not unique, and suppose further that $\hat{y}$ is
maximal among points for where the supremum is attained other than
$y$. 

Take any $q \in \mathbb{Q}$ with $\hat{y} < q < y$. The remainder of
the argument is a time reversed analogue of the argument used in the
proof of Theorem~\ref{thm:Aclosed} with a slightly expanded defnition
of $E^q$. Let $\bar{C}_q : [q, \infty) \ra \reals$ be
the concave majorant of $\left( \psi(x) - \frac{1}{2} x^2 \right)_{ x
  \geq q}$ and let $\bar{c}_q$ be its right continuous derivative,
which will agree with $\bar{c}$ on the set $(y, \infty)$. Define
\[
E_q : = \{ x \geq q \, : \, \bar{c}_q(x-) \leq -x \, , \, \bar{c}_q(x)
\geq x \} \, . 
\]
Since $x \in E_q$ implies that at least one of $\limsup_{h \downarrow 0} h^{-1} (
\psi_0(x+h) - \psi_0(x) )$ or $\limsup_{h \downarrow 0} h^{-1} (
\psi_0(x-h) - \psi_0(x-) )$ is finite, \eqref{eq:hypB} and Fubini
imply that $E_q$ has measure zero almost surely. Also, $a(y) = y$ it
follows that $y \in E_q$.

For every $x \in E_q$ define
\[
T_q(x) := \inf \left\{ z \ge 0 \, : \, \psi_0(q-z) \vee
\psi_0((q-z)-) - \mbox{$\frac{1}{2}$} (q-z)^2 \geq \psi_0(x) \vee
\psi_0(x-) - x \left( (q-x) - z \right) \right\} \, ,
\]
and note that almost surely $T_q(x) > 0$ for every $x \in
E_q$ since by \eqref{eq:hypB} $\bar{C}(q) > \psi_0(q) \vee
\psi_0(q-) - \frac{1}{2} q^2$ a.s. Also, by definition $T_q(y) = \hat{y}$.

By Assumption~\textbf{B}(II) (its time reversed version -- see
Remark~\ref{rem:assB}(ii)) and the fact that $E_q$ has measure zero
a.s. it follows that a.s.
\begin{equation}
\label{eq:yhat2}
\psi_0(T_q(x)) - \mbox{$\frac{1}{2}$} (q-T_q(x))^2 \geq \psi_0(x) \vee \psi_0(x-) - y \left( (q-x) - T_q(x) \right)
\end{equation}
for every $x \in E_q$ such that $T_q(x) < \infty$. But $y \in
E_q$ and $T_q(y) = \hat{y}$ a.s. hence \eqref{eq:yhat} and
\eqref{eq:yhat2} would imply that $\hat{y} = - \infty$, and thus
$\hat{y}$ cannot exist as assumed. 
\endproof

\section{Proof of Theorem~\ref{thm:regenerative}}
\label{sec:regenproof}

\subsection{Facts relating to the first non-negative element of $\mathcal{A}_0$}
\label{sec:facts}

In this section, we prove some results relating to the first
non-negative element of $\mathcal{A}_0$ when $\psi_0$ is a non-random
c\`adl\`ag function satisfying $\lim_{|x| \ra \infty} x^{-2} \psi_0(x) = 0$.
Define
\[
\begin{split}
\mathbf{t} : & \! \! = \inf \{ y \geq 0 \, : \, a(y) = y \} \\
& \! \! = \inf \left\{ y \geq 0 \, : \,
\arg \sup \left\{ \psi_0(x) - \mbox{$\frac{1}{2}$} (x-y)^2 : x \in \reals \right\} 
= y \right\} , \\
& \! \! = \inf \left\{ y \geq 0 \, : \,
\psi_0(y-x) \vee \psi_0((y-x)-) - \psi_0(y) \vee \psi_0(y-) \leq \mbox{$\frac{1}{2}$} x^2
\text{ for all } x > 0 \text{ and }  \color{white} \right\}   \\
& \quad \quad  \quad \quad  \quad \, \, \color{white} \left\{ \color{black} \psi_0(y+x) \vee \psi_0((y+x)-)  - \psi_0(y) \vee \psi_0(y-) < \mbox{$\frac{1}{2}$} x^2
\text{ for all } x > 0  \right\} . \\
\end{split}
\]  
The last equality is becuase of the convention that if the $\arg \sup$
above is not unique we take it to be the supremum over all suitable
arguments. Define further
\[
\begin{split}
\mathbf{r} & := \inf \left\{ y \geq 0 \, : \, \psi_0(y-x) \vee \psi_0((y-x)-) - \psi_0(y) \vee \psi_0(y-) \leq \mbox{$\frac{1}{2}$} x^2
\text{ for all } x > 0 \right\} ,\\
\mathbf{s} & := \inf \left\{ y \geq \mathbf{r} \, : \, \psi_0(y+x) \vee \psi_0((y+x)-)  - \psi_0(y) \vee \psi_0(y-) < \mbox{$\frac{1}{2}$} x^2
\text{ for all } x > 0 \right\} . \\
\end{split}
\]
Note that $0 \leq \mathbf{r} \leq \mathbf{s} \leq \mathbf{t}$.

\begin{lem}
\label{lem:sinf}
Let $\psi_0$ be any c\`adl\`ag function with $\lim_{|x| \ra
  \infty} x^{-2} \psi_0(x) = 0$.
Then the infimum in the definition of $\mathbf{s}$ is achieved, that
is,
\begin{equation}
\label{eq:sinf}
\psi_0(\mathbf{s}+x) \vee \psi_0((\mathbf{s}+x)-)  -
\psi_0(\mathbf{s}) \vee \psi_0(\mathbf{s}-) < \mbox{$\frac{1}{2}$} x^2
\text{ for all } x > 0 \, .
\end{equation}
\end{lem}
\proof
Suppose that \eqref{eq:sinf} did not hold. Then by the definition of
$\mathbf{s}$ there would exist a strictly decreasing sequence $\{
\mathbf{s}_n   \}_{n \ge 0}$ such that $\lim_n \mathbf{s}_n =
\mathbf{s}$ and
\begin{equation}
\label{eq:xn}
\psi_0(\mathbf{s}_n+x) \vee \psi_0((\mathbf{s}_n+x)-)  - \psi_0(\mathbf{s}_n) \vee \psi_0(\mathbf{s}_n-) <
\mbox{$\frac{1}{2}$} x^2 \text{ for all } x > 0 
\end{equation}
for every $n \ge 0$. 

For $n \ge 1$, \eqref{eq:xn} with $x = \mathbf{s}_{n-1} -
\mathbf{s}_n$ gives
\[
\psi_0(\mathbf{s}_n) \vee \psi_0(\mathbf{s}_n-) 
> 
\psi_0(\mathbf{s}_{n-1}) \vee \psi_0(\mathbf{s}_{n-1}-) 
- \mbox{$\frac{1}{2}$}(\mathbf{s}_{n-1} - \mathbf{s}_n)^2 
\]
and thus since $\sum_{m=1}^n (\mathbf{s}_{m-1} - \mathbf{s}_m)^2 \leq
(\mathbf{s}_0 - \mathbf{s}_n)^2$ we have
\[
\psi_0(\mathbf{s}_n) \vee \psi_0(\mathbf{s}_n-) 
> 
\psi_0(\mathbf{s}_0) \vee \psi_0(\mathbf{s}_0-) 
- \mbox{$\frac{1}{2}$}(\mathbf{s}_0 - \mathbf{s}_n)^2 .
\]
By right continuity of $\psi_0(\cdot)$ at $\mathbf{s}$, recalling that
$\lim_n \mathbf{s}_n = \mathbf{s}$ we may take the limit as $n \ra
\infty$ to get that 
\begin{equation}
\label{eq:sands0}
\psi_0(\mathbf{s}) 
\geq
\psi_0(\mathbf{s}_0) \vee \psi_0(\mathbf{s}_0-) -  \mbox{$\frac{1}{2}$}(\mathbf{s}_0 - \mathbf{s})^2 .
\end{equation}

Now, since we have assumed that \eqref{eq:sinf} does not hold, there
exists $x^*>0$ such that
\[
\psi_0(\mathbf{s}+x^*) \vee \psi_0((\mathbf{s}+x^*)-)  -
\psi_0(\mathbf{s}) \vee \psi_0(\mathbf{s}-) \geq \mbox{$\frac{1}{2}$}
(x^*)^2 \, ,
\]
and moreover wihtout loss of generality we can assume that
$\mathbf{s}_0$ is such that $\mathbf{s}_0 < \mathbf{s} + x^*$.
But then starting from \eqref{eq:sands0} we get
\[
\begin{split}
\psi_0(\mathbf{s}_0) \vee \psi_0(\mathbf{s}_0-)
& \leq
\psi_0(\mathbf{s}) \vee \psi_0(\mathbf{s}-) + \mbox{$\frac{1}{2}$}(\mathbf{s}_0-\mathbf{s})^2 \\
& \leq
\psi_0(\mathbf{s}) \vee \psi_0(\mathbf{s}-) + \mbox{$\frac{1}{2}$}(x^*)^2
- \mbox{$\frac{1}{2}$} ((\mathbf{s}+x^*)-\mathbf{s}_0)^2  \\
& \leq \psi_0(\mathbf{s}+x^*) \vee \psi_0((\mathbf{s}+x^*)-)
- \mbox{$\frac{1}{2}$} ((\mathbf{s}+x^*)-\mathbf{s}_0)^2 , \\
\end{split}
\]
which contradicts \eqref{eq:xn} with $n=0$ and $x = \mathbf{s}+x^*-\mathbf{s}_0$.
\endproof


\begin{lem}
\label{lem:s=t}
Let $\psi_0$ be any c\`adl\`ag function with $\lim_{|x| \ra \infty} x^{-2} \psi_0(x) = 0$. Then $\mathbf{s} = \mathbf{t}$.
\end{lem}
\proof
Recall that $0 \leq \mathbf{r} \leq \mathbf{s} \leq \mathbf{t}$. We
will show that at $\mathbf{s}$ the conditions of $\mathbf{r}$ are
still satisfied, i.e.
\begin{equation}
\label{eq:sliker}
\psi_0(\mathbf{s}-x) \vee \psi_0((\mathbf{s}-x)-) - \psi_0(\mathbf{s})
\vee \psi_0(\mathbf{s}-) \leq \mbox{$\frac{1}{2}$} x^2 \text{ for all } x > 0
\end{equation}
which combined with~\eqref{eq:sinf} implies that $\mathbf{s} \geq
\mathbf{t}$ and hence $\mathbf{s} = \mathbf{t}$.

Suppose first that $\mathbf{r} = \mathbf{s}$, then
clearly~\eqref{eq:sliker} is satisfied and hence $\mathbf{s}=\mathbf{t}$.
Assume therefore that $\mathbf{r} < \mathbf{s}$.
We will begin by showing that \eqref{eq:sliker} holds for all $0 < x
\leq \mathbf{s}-\mathbf{r}$. It suffices to show that if we define
\begin{equation}
\label{eq:s=t_tau}
\tau := \arg \sup \left\{ \psi_0(\mathbf{s}-y) \vee \psi_0((\mathbf{s}-y)-) - \mbox{$\frac{1}{2}$}
y^2 : 0 \leq y \le \mathbf{s}-\mathbf{r} \right\} 
\end{equation}
then we must have $\tau = 0$. Well, \eqref{eq:s=t_tau} implies that
\[
\psi_0(\mathbf{s}-\tau) \vee \psi_0((\mathbf{s}-\tau)-) - \mbox{$\frac{1}{2}$} \tau^2 
\ge
\psi_0(\mathbf{s}-y) \vee \psi_0((\mathbf{s}-y)-) - \mbox{$\frac{1}{2}$} y^2 
\]
for all $0 \leq y \leq \tau$. 
Making the change of variables $y = \tau - x$, we see that
\[
\psi_0(\mathbf{s}-\tau + x) \vee \psi_0((\mathbf{s}-\tau + x)-) 
- \psi_0(\mathbf{s}-\tau) \vee \psi_0((\mathbf{s}-\tau)-) 
\leq 
\mbox{$\frac{1}{2}$} x^2 - x \tau
\]
for all $0 \le x \le \tau$. Suppose that $\tau > 0$, so that
\[
\psi_0(\mathbf{s}-\tau + x) \vee \psi_0((\mathbf{s}-\tau + x)-) 
- \psi_0(\mathbf{s}-\tau) \vee \psi_0((\mathbf{s}-\tau)-) 
<
\mbox{$\frac{1}{2}$} x^2
\]
for all $0 < x \le \tau$. Combined with \eqref{eq:sinf} this would
imply that
\[
\psi_0(\mathbf{s}-\tau) \vee \psi_0((\mathbf{s}-\tau)-) + \mbox{$\frac{1}{2}$}
x^2 > \psi((\mathbf{s}-\tau)+x) \text{ for all } x > 0 \, .
\]
But then since $\mathbf{s}-\tau \ge \mathbf{r}$, the definition of $\mathbf{s}$
would then imply that $\mathbf{s} \le \mathbf{s}-\tau < \mathbf{s}$, a clear contradiction. Hence
$\tau = 0$ as required.

It remains to show that \eqref{eq:sliker} holds for all $x > \mathbf{s}-\mathbf{r}$.
Applying \eqref{eq:sliker} at $x = \mathbf{s}-\mathbf{r}$ we see that
\[
\psi_0(\mathbf{r}) \vee \psi_0(\mathbf{r}-) - \psi_0(\mathbf{s}) \vee \psi_0(\mathbf{s}-) 
\leq 
\mbox{$\frac{1}{2}$} (\mathbf{s}-\mathbf{r})^2. 
\]
From the definition of $\mathbf{r}$, 
\[
\psi_0(\mathbf{r}-y) \vee \psi_0((\mathbf{r}-y)-) - \psi_0(\mathbf{r}) \vee \psi_0(\mathbf{r}-) \leq \mbox{$\frac{1}{2}$} y^2
\]
for all $y>0$, and hence
\[
\psi_0(\mathbf{r}-y) \vee \psi_0((\mathbf{r}-y)-) - \psi_0(\mathbf{s}) \vee \psi_0(\mathbf{s}-) 
\leq
\mbox{$\frac{1}{2}$} y^2 + \mbox{$\frac{1}{2}$} (\mathbf{s}-\mathbf{r})^2
<
\mbox{$\frac{1}{2}$} ((\mathbf{s}-\mathbf{r})+y)^2
\]
for all $y>0$. Applying the change of variables $x =
(\mathbf{s}-\mathbf{r}) + y$ shows that \eqref{eq:sliker} holds for
all $x > \mathbf{s}-\mathbf{r}$ and hence completes the proof.
\endproof

Define $\mathbf{r}_0 := \mathbf{r}$, and for $k \ge 0$, define
\[
\mathbf{r}_{k+1} := \mathbf{r}_k
+ \arg \sup
\left\{ 
\psi_0(\mathbf{r}_k+x) \vee \psi_0((\mathbf{r}_k+x)-) - \mbox{$\frac{1}{2}$} x^2
: x \ge 0 \right\}  \, ,
\]
where if the $\arg \sup$ above is not unique we take it to be the supremum
over all suitable arguments.

\begin{lem}
\label{lem:r_ktos}
Let $\psi_0$ be any c\`adl\`ag function with $\lim_{|x| \ra \infty}
x^{-2} \psi_0(x) = 0$. Then $\mathbf{r}_k \ra \mathbf{t}$.
\end{lem}
\proof
Note first that $\mathbf{r}^* := \lim_k \mathbf{r}_k$ exists
 since $\mathbf{r}_k$ is an increasing sequence.
If there were a $k \ge 0$ such that $\mathbf{r}_k = \mathbf{t}$, then
necessarily $\mathbf{r}_j = \mathbf{t}$ for all $j \ge k$, thus we henceforth assume
there is no such $k$.  

Suppose that there exists a $k \ge 0$ such that $\mathbf{r}_k < \mathbf{t} < \mathbf{r}_{k+1}$,
then 
\begin{equation}
\label{eq:rlim1}
\psi_0(\mathbf{t}) \vee \psi_0(\mathbf{t}-) -
\mbox{$\frac{1}{2}$} (\mathbf{t}-\mathbf{r}_k)^2
 \leq
\psi_0(\mathbf{r}_{k+1}) \vee \psi_0(\mathbf{r}_{k+1}-) -
\mbox{$\frac{1}{2}$} (\mathbf{r}_{k+1}-\mathbf{r}_k)^2. 
\end{equation}
From the definition of $\mathbf{t}$ it follows that
\[
\psi_0(\mathbf{r}_{k+1}) \vee \psi_0(\mathbf{r}_{k+1}-) - \psi_0(\mathbf{t}) \vee \psi_0(\mathbf{t}-) -
\mbox{$\frac{1}{2}$} (\mathbf{r}_{k+1}-\mathbf{t})^2 < 0 .
\]
Thus if equality held in \eqref{eq:rlim1} it would be the case
that
\[
(\mathbf{r}_{k+1}-\mathbf{r}_k)^2 < (\mathbf{r}_{k+1}-\mathbf{t})^2 + (\mathbf{t}-\mathbf{r}_k)^2 
= (\mathbf{r}_{k+1}-\mathbf{r}_k)^2 - 2(\mathbf{r}_{k+1}-\mathbf{t})(\mathbf{t}-\mathbf{r}_k),
\]
and hence the inequality in \eqref{eq:rlim1} must be strict.
\eqref{eq:rlim1} then implies that
\[
\begin{split}
& \psi_0(\mathbf{r}_{k+1}) \vee \psi_0(\mathbf{r}_{k+1}-) - \psi_0(\mathbf{t}) \vee \psi_0(\mathbf{t}-) 
- \mbox{$\frac{1}{2}$} (\mathbf{r}_{k+1}-\mathbf{t})^2 \\
& \qquad \qquad \qquad \qquad \qquad >
(\mathbf{r}_{k+1}-\mathbf{r}_k)^2 - (\mathbf{r}_{k+1}-\mathbf{t})^2 - (\mathbf{t}-\mathbf{r}_k)^2 > 0 ,\\
\end{split}
\]
which contradicts the definition of $\mathbf{t}$, and hence there is no $k$
such that $\mathbf{r}_k < \mathbf{t} < \mathbf{r}_{k+1}$. Thus
$\mathbf{r}^* \leq \mathbf{t}$. 

Suppose $\mathbf{r}^* < \mathbf{t}$, then $\mathbf{r}^* < \mathbf{s}$ by Lemma~\ref{lem:s=t},
and hence there exists $r_+>0$ such that
\[
\psi_0(\mathbf{r}^*+r_+) - \psi_0(\mathbf{r}^*) \vee \psi_0(\mathbf{r}^*-) -
\mbox{$\frac{1}{2}$} r_+^2 > 0.
\]
Let $r_- > 0$ be such that 
\[
\mbox{$\frac{1}{2}$} (r_+ + r_-)^2 = \psi_0(\mathbf{r}^*+r_+) - \psi_0(\mathbf{r}^*) \vee \psi_0(\mathbf{r}^*-) .
\]
Then for all $k$ large enough such that $\mathbf{r}_k > \mathbf{r}^* - r_-$ we have
\[
\mbox{$\frac{1}{2}$} ((\mathbf{r}^*+r_+) - \mathbf{r}_k)^2 < \psi_0(\mathbf{r}^*+r_+) - \psi_0(\mathbf{r}^*) \vee \psi_0(\mathbf{r}^*-)
\]
and hence
\[
\mathbf{r}_{k+1} > \mathbf{r}_k
+ ((\mathbf{r}^*+r_+) -\mathbf{r}_k) = \mathbf{r}^* + r_+ > \mathbf{r}^*,
\]
which is clearly a contradiction. Thus we can conclude that $\mathbf{r}^*=\mathbf{t}$.
\endproof

\subsection{Randomized coterminal times relating first non-negative element of $\mathcal{A}_0$}
\label{sec:prelims}

In this section we will use the notation of Definition~\ref{def:rct}
when checking if a given random time is a randomized coterminal time. 

\begin{lem}
\label{lem:r_k_rct}
Let $\psi_0$ be a real valued strong Markov process. Define a sequence
of random times by $R_0 = 0$, and for $k \ge 0$,
\[
\begin{split}
R_{k+1} &  := R_k
+ \arg \sup
\left\{ 
\psi_0(R_k+x) \vee \psi_0((R_k+x)-) - \psi_0(R_k) \vee
\psi_0(R_k) - \mbox{$\frac{1}{2}$} x^2
: x \ge 0 \right\} \, \\
& = R_k
+ \arg \sup
\left\{ 
\psi_0(R_k+x) \vee \psi_0((R_k+x)-) - \mbox{$\frac{1}{2}$} x^2
: x \ge 0 \right\} \, , \\
\end{split}
\]
where if the $\arg \sup$ above is not unique we take it to be the
supremum over all suitable arguments.
Define $\phi$ to be the process
$(\phi(x))_{x \ge 0}$, with
\[
\phi(x) = (\phi_1(x),\phi_2(x)) := (x, \psi_0(x))
\]
for all $x \ge 0$. Then $R_k$ is a randomized coterminal time for
$\phi$ for each $k \ge 1$.
\end{lem}
\proof
Let $A = \reals^3$, $\mathfrak{U} = \mathcal{B}(\reals^3)$ and let
$Z = (R_{k-1}, R_k, \psi_0(R_k) \vee \psi_0(R_k-))$. Let $R_0^{(x)} =
0$ and for $k \ge 0$, if $R_k^{(x)} < x$ then let
\[
\begin{split}
R_{k+1}^{(x)} & := R_k^{(x)} -  
+ \arg \sup
\left\{ 
\psi_0(R_k^{(x)}+y) \vee \psi_0((R_k^{(x)}+y)-) - \mbox{$\frac{1}{2}$} y^2
: 0 \le y \le x - R_k^{(x)} \right\} , \\
\end{split}
\]
but if $R_k^{(x)} = x$ then let $R_{k+1}^{(x)}= x$.
Let $Z_x = (R_{k-1}^{(x)}, R_k^{(x)}, \psi_0(R_k^{(x)}) \vee \psi_0(R_k^{(x)}-))$, 
so that $Z_x$ is an $\mathcal{F}_x$-measurable $A$-valued
random variable as required. Finally, recalling that
$(\phi_1(x),\phi_2(x)) = (x, \psi_0(x))$, define the family of
terminal times $\{ T_a \}_{a \in A}$ by
\[
T_{(a_1,a_2,a_3)} := 
\inf \{ x > 0 : \phi_2(x) \vee \phi_2(x-) 
- a_3 + \mbox{$\frac{1}{2}$} (a_1-a_2)^2
\geq \mbox{$\frac{1}{2}$} ( \phi_1(x) - a_1)^2 \}.
\]
(I) and (II) follow once we define $B(y,x) := \{ y \le R_k^{(x)} < x \}$.
\endproof

\begin{lem}
\label{lem:first_forward_parabola}
Let $\psi_0$ be a real valued strong Markov process and define
\[
F := \inf \left\{ x \geq 0 \, : \, \psi_0(x+s) \vee \psi_0((x+s)-) - \psi_0(x) \vee
  \psi_0(x-) < \mbox{$\frac{1}{2}$}  s^2 \text{ for all } s > 0 \right\} .
\] 
Define $\phi$ to be the process $(\phi(x))_{x \ge 0}$, with 
\[
\phi(x) = (\phi_1(x),\phi_2(x)) := (x, \psi_0(x))
\]
for all $x \ge 0$. Then $F$ is a randomized coterminal time for $\phi$.
\end{lem}
\proof
Let $A = \reals^2$, $\mathfrak{U} = \mathcal{B}(\reals^2)$,
$Z = (F, \psi_0(F) \vee \psi_0(F-))$ and $Z_x = (F_x, \psi_0(F_x) \vee
\psi_0(F_x-))$, where
\[
F_x := \inf \{ 0 \leq y \leq x \, : \, 
\psi_0(y+s) \vee \psi_0((y+s)-)  - \psi_0(y) \vee \psi_0(y-) 
\leq \mbox{$\frac{1}{2}$}  s^2 \text{ for all } 0 < s \le x-y \} .
\]
It follows that $Z_x$ is an $\mathcal{F}_x$-measurable $A$-valued
random variable. Finally, recalling that
$(\phi_1(x),\phi_2(x)) = (x, \psi_0(x))$, define the family of
terminal times $\{ T_a \}_{a \in A}$ by
\[
T_{(a_1,a_2)} := 
\inf \{ x > 0 : \phi_2(x) \vee \phi_2(x-) - a_2
\geq \mbox{$\frac{1}{2}$} (\phi_1(x)-a_1)^2 \}.
\]

By definition,
\[
\psi_0(F+s) \vee \psi_0((F+s)-) - \psi_0(F) \vee
  \psi_0(F-) < \mbox{$\frac{1}{2}$}  s^2 
\]
for all $s > 0$, and
\[
\psi_0(F_x+s) \vee \psi_0((F_x+s)-) - \psi_0(F_x) \vee \psi_0(F_x-) < \mbox{$\frac{1}{2}$}  s^2
\]
for all $0 < s \leq x-F_x$. In particular, if $F \leq x$, then
\[
\psi_0(F) \vee \psi_0(F-)  - \psi_0(F_x) \vee \psi_0(F_x-) <
\mbox{$\frac{1}{2}$}  (F-F_x)^2 .
\]
Hence we see that on the set $\{ F \leq x \}$, 
\[
\psi_0(F_x+s) \vee \psi_0((F_x+s)-) - \psi_0(F_x) \vee \psi_0(F_x-) < \mbox{$\frac{1}{2}$}  s^2
\]
for all $s > 0$, which implies that $F \leq F_x$. However, $F \geq
F_x$ by definition, and therefore $F_x = F$ on the set $\{ F \leq x
\}$. Thus (I) is satisfied.

If we define $B(y,x) := \{ y \le F_x < x \}$, then clearly
\[
\{ y \leq F < x \}
=
B(y,x) \cap \{ T_{Z_x(\omega)} (\theta_x \omega) = +\infty \} 
=
B(y,x) \cap \{ T_{Z(\omega)} (\theta_x \omega) = +\infty \} ,
\]
and hence (II) is satisfied.
\endproof

\begin{cor}
\label{cor:first_forward_parabola}
Let $\psi_0$ be a L\'evy process and define $F$ as in
Lemma~\ref{lem:first_forward_parabola}. Suppose that
$\psi_0$ is continuous at $F$. Then for any $(x_1,\ldots,x_n)$ with
$x_i > 0$ for $i = 1,\ldots,n$, the joint law of
$(\psi_0(F+x_i)-\psi_0(F))_{i = \{ 1, \ldots , n \} }$ depends only on
  $(x_1,\ldots,x_n)$.
\end{cor}
\proof
From Theorem~\ref{thm:rct_markov} we know that the joint
law of $(\psi_0(F+x_i))_{i = \{ 1, \ldots , n \} }$ depends only on
$(x_1,\ldots,x_n)$ and $Z = (F, \psi_0(F))$. Moreover we can think of
the post $F$ process $(\psi_0(F+x))_{x \ge 0}$ as the original process
started at $\psi_0(F)$ but conditioned to remain below a half parabola
with its minimum at $\psi_0(F) \vee \psi_0(F-) = \psi_0(F)$. Then by
the spatial homogeneity of L\'evy processes, the joint law of
$(\psi_0(F+x_i)-\psi_0(F))_{i = \{ 1, \ldots , n \} }$ cannot depend
on $\psi_0(F)$, and by the temporal homogeneity of L\'evy processes it
cannot depend on $F$ either. Thus the joint law of $(\psi_0(F+x_i) -
\psi_0(F))_{i = \{ 1, \ldots , n \} }$ can depend only on $(x_1,\ldots,x_n)$.
\endproof

\subsection{Proof of Theorem~\ref{thm:regenerative}}
\label{sec:regenproofsubsection}

\proof
\textit{(Theorem~\ref{thm:regenerative})}
  Recall from the statement of the theorem that $T := \inf \{ x \geq 0
  \, : \, x \in \mathcal{A}_0 \} $ and hence $T = \inf \{ x \geq 0 \,
  : \, a(x) = x \}$.  From Corollary~\ref{cor:bounded} in the bounded
  variation case or Theorem~\ref{thm:abrupt} in the abrupt case, we
  know that the set $\{ x \in \reals \, : \, a(x) = x \}$ is closed
  a.s. and hence $a(T) = T$ a.s.

  If $\psi_0$ has paths of unbounded variation, then
  Lemma~\ref{lem:unbounded_continuous} then implies that $\psi_0$ is
  continuous at $T$ a.s.  If $\psi_0$ has paths of bounded variation
  Theorem~\ref{thm:zeros_are_regular} implies that $T$ is a Lagrangian
  regular point a.s. and then Proposition~\ref{prop:bounded}(ii) implies
  that $\psi_0$ is continuous at $T$ a.s.

Define two further random variables $R$ and $S$ by
\[
\begin{split}
R & := \inf \left\{ y \geq 0 \, : \, \psi_0(y-x) \vee \psi_0((y-x)-) - \psi_0(y) \vee \psi_0(y-) \leq \mbox{$\frac{1}{2}$} x^2
\text{ for all } x > 0 \right\} ,\\
S & := \inf \left\{ y \geq R \, : \, \psi_0(y+x) \vee \psi_0((y+x)-)  - \psi_0(y) \vee \psi_0(y-) < \mbox{$\frac{1}{2}$} x^2
\text{ for all } x > 0 \right\} . \\
\end{split}
\]
Note that $0 \leq R \leq S \leq T$.
Note also that by the strong Markov property
applied at the stopping time $R$, it follows that $S-R$ has the same
law as $F$ in Lemma~\ref{lem:first_forward_parabola}.

Lemma~\ref{lem:s=t} tells us that $S=T$ a.s., and thus $T = R + (S-R)$
a.s. Since $R$ is a stopping time, $(\psi_0(R+x)-\psi_0(R))_{x \ge 0}$
is independent of $(\psi_0(R-x))_{x \ge 0}$ and has the same law as
$(\psi_0(x))_{x \ge 0}$. Since $S-R$ has the same law as $F$ in
Lemma~\ref{lem:first_forward_parabola}, we only need to show that the
process $(\psi_0(F+x)-\psi_0(F))_{x \ge 0}$ is independent of
$(\psi_0(x))_{0 \le x \le F}$ when $\psi_0$ is a.s. continuous at
$F$, and when we can further assume that
\begin{equation}
\label{eq:convenient}
\psi_0(x) - \mbox{$\frac{1}{2}$} x^2 \leq 0 \text{  for all $x \leq 0$.}
\end{equation}

By continuity of $\psi_0$ at $F$ we only need to show
that $(\psi_0(F+x)-\psi_0(F))_{x \ge 0}$ is independent of
$(\psi_0(x))_{0 \le x < F}$. Moreover,
Corollary~\ref{cor:first_forward_parabola} implies that the law of
$(\psi_0(F+x)-\psi_0(F))_{x \ge 0}$ cannot depend on $F$ or
$\psi_0(F)$, hence it is enough to show that 
\begin{equation}
\label{eq:show_this}
\begin{split}
& \text{$(\psi_0(F+x)-\psi_0(F))_{x \ge 0}$ is independent of
$(\psi_0(x))_{0 \le x < F}$} , \\
& \qquad \qquad \qquad \qquad \qquad
\qquad \qquad \qquad
\text{   conditionally given $F$ and $\psi_0(F)$.} \\
\end{split}
\end{equation}

Suppose first that $\psi_0$ has paths of unbounded variation and is
abrupt. From Corollary~\ref{cor:local_max} $\psi_0$
must have a local maximum at $F$, and from
Lemma~\ref{lem:first_forward_parabola} we know that $F$ is a 
randomized coterminal time for the process $(x, \psi_0(x))_{x \ge 0}$, hence
by Proposition~\ref{prop:local_max_and_jump_rct} it follows that 
$(F+x,\psi_0(F+x)-\psi_0(F))_{x \ge 0}$
is independent of $(x,\psi_0(x))_{0 \le x \le F}$
conditionally given $(F, \psi_0(F))$. Hence we
have~\eqref{eq:show_this}.

Now suppose that $\psi_0$ has paths of bounded variation. 
Define a sequence of random times by $R_0 = 0$, and
\[
R_{k+1} := R_k
+ \arg \sup
\left\{ 
\psi_0(R_k+x) \vee \psi_0((R_k+x)-) - \mbox{$\frac{1}{2}$} x^2
: x \ge 0 \right\} 
\]
for $k \ge 0$. 

From Lemma~\ref{lem:r_ktos} and the fact that $S-R$ has
the same law as $F$, we have that $R_k \ra F$ a.s. Suppose we have shown that
for each $k \ge 1$, the process $(\psi_0(R_k+x))_{ x  \ge 0 }$ is
independent of $(\psi_0(x))_{0 \le  x  \le R_k }$ conditionally given
$R_k$ and $\psi_0(R_k)$. If $R_k = R_{k+1}$ for some $k$, then $R_k =
T$ and we are done. Thus assume that $R_k < R_{k+1}$ for every $k$. 
We have that $R_k \ra F$ a.s., and the
a.s. continuity of $\psi_0$ at $F$ implies that $\psi_0(R_k) \ra
\psi_0(F)$. Thus the process $(\psi_0(F+x))_{ x  \ge 0 }$ is
independent of $(\psi_0(x))_{0 \le  x  < R_k }$ conditionally given
$F$ and $\psi_0(F)$, and \eqref{eq:show_this} follows.

It remains to show that for each $k \ge 1$, the process $(\psi_0(R_k+x))_{ x  \ge 0 }$ is
independent of $(\psi_0(x))_{0 \le  x  \le R_k }$ conditionally given
$R_k$ and $\psi_0(R_k)$. 
Note that under~\eqref{eq:convenient}, $a(R_k) = R_{k+1}$ for every $k
\ge 0$. Since we have assumed that $R_k < R_{k+1}$ for every $k \ge
0$, it follows from Proposition~\ref{prop:l-rey} that $R_{k+1}$ is a
positive jump time of $\psi_0$ for every $k \ge 0$. From
Lemma~\ref{lem:r_k_rct} we know that $R_k$ is a randomized coterminal
time for the process $(x, \psi_0(x))_{x \ge 0}$, hence
by Proposition~\ref{prop:local_max_and_jump_rct} it follows that 
$(R_k+x,\psi_0(R_k+x)-\psi_0(R_k))_{x \ge 0}$
is independent of $(x,\psi_0(x))_{0 \le x \le R_k}$
conditionally given $(R_k, \psi_0(R_k))$. 
\endproof

\begin{remark}
For processes with bounded variation satisfying Hypotheses~\textbf{A}
and \textbf{B}, Giraud's proof of the regenerativity of the set of
Lagrangian regular points
\cite[Theorem 2]{giraud_burgers} when $\psi_0$ is a stable L\'evy
process with stability index $\alpha \in (1/2,1)$ could also be used
to prove Theorem~\ref{thm:regenerative}. Hypothesis~\textbf{B} ensures
that equation (7) of \cite{giraud_burgers} holds appropriately, and
Theorem~\ref{thm:zeros_are_regular} ensures that the first sentence of
Lemma 4 of \cite{giraud_burgers} is true. Those are the only two
results needed in that proof.
\end{remark}

\section{Acknowledgements}
\label{sec:thanks}

The author would like to thank Steve Evans for generating many of the
ideas used in this paper during the joint writing of~\cite{lipz} and
Jim Pitman for suggesting reading that lead to the author working on
this topic. In addition, the author is very grateful for the detailed
comments and helpful suggestions of the anonymous referees.

\bibliographystyle{plain}
\bibliography{parabolic_arxiv_revised.bib}

\end{document}